\documentclass[11pt]{article}
\usepackage[english]{babel}
\usepackage[utf8]{inputenc}
\usepackage{bm}
\usepackage{hyperref} 
\usepackage{amsmath,amsthm,amsfonts}
\usepackage{graphicx}
\usepackage{cite}
\usepackage{xcolor}
\usepackage{multirow}
\usepackage{subcaption}
\usepackage{fancyhdr}  
\usepackage{geometry}  
\usepackage{caption}
\usepackage{makecell} 
\usepackage{array}
\usepackage{float} 
\usepackage{enumitem} 
\usepackage{titlesec}
\usepackage{booktabs} 
\usepackage{lineno}

\usepackage{url} 
\graphicspath{ {./figs/} }

\usepackage{todonotes}

\newcommand{\mwang}[1]{{\color{black} #1}}
\newcommand{\yw}[1]{{\color{black} #1}}
\date{}

\newcommand{\mytitle}{Numerical Methods for Solving Nonlinearly Coupled Poisson Equations in Dual-Continuum Modeled Porous Electrodes}

\title{\mytitle}

\author{
Yuhe Wang\thanks{Institute for Scientific Computation, Texas A\&M University, College Station, TX 77843, USA. \texttt{yuhe.wang@tamu.edu}} \and
Min Wang\thanks{Department of Mathematics, University of Houston, Houston, TX 77204, USA. \texttt{mwang55@central.uh.edu}} \thanks{Corresponding author} \and
Zhihang Xu\thanks{Department of Mathematics, University of Houston, Houston, TX 77204, USA. \texttt{zxu29@central.uh.edu}}
}

\date{}

\begin{document}

\maketitle

\thispagestyle{empty}

\pagestyle{fancy}  

\begin{abstract}
Porous electrodes are widely used in electrochemical systems, where accurately determining electric potentials, particularly overpotentials, is essential for understanding electrode behavior. At the macroscopic scale, we model porous electrodes using a dual-continuum formulation, in which the porous solid phase and the liquid electrolyte are represented as spatially superimposed continua. Determining the potential distributions then requires solving two Poisson equations that are nonlinearly coupled through Butler--Volmer kinetics under galvanostatic and potentiostatic operating modes. Under galvanostatic operation, these equations form an underconstrained singular system due to all-Neumann boundary conditions, posing numerical challenges. This paper systematically presents the numerical methods for solving the resulting nonlinearly coupled Poisson system in a dual-continuum framework, with a particular focus on galvanostatic solutions. We establish solution uniqueness in terms of the potential difference between the electrode and electrolyte (or overpotential), as well as the individual potentials up to a shared constant shift. To resolve this nonuniqueness, we introduce three numerical approaches: (1) Lagrange constrained method (LCM), (2) Dirichlet substitution method (DSM), and (3) global constraining method (GCM), which enables solving the overpotential without imposing an explicit system reference potential. Additionally, we provide both decoupled and fully coupled nonlinear solution strategies and evaluate their computational performance for homogeneous and heterogeneous conductivity fields. The presented modeling and numerical framework is general and can be applied to related underconstrained nonlinear systems.
\vspace{1em}  

\end{abstract}

\section{Introduction}
\vspace{-1.8ex}
Porous electrodes play a crucial role in many electrochemical energy storage devices \cite{SANTHANAGOPALAN2009110,newman2021electrochemical}. Their \mwang{large} surface area and interconnected pore structures facilitate efficient electrochemical reactions and mass transport, both of which are critical to device performance. The efficiency of porous electrodes arises from the intricate interplay among ionic transport in the electrolyte, electronic conduction in the solid electrode phase, and electrochemical reactions occurring at the electrode-electrolyte interfaces within the porous structure \cite{alkire1975flow, pan2022porous, hamed2021limitation, qu2014fundamental}. Understanding and optimizing these coupled processes are essential for advancing the design and performance of electrochemical energy storage devices, including batteries, flow cells, and supercapacitors \cite{smith2017multiphase, zhang2020understanding}. A key aspect of this is accurately modeling the electric potential fields in both the solid electrode and liquid electrolyte phases, as these fields dictate the localized distribution of electrochemical reaction rates and ultimately influence overall device efficiency
\cite{newman1975porous,AMINI2022466}. 

The continuum approach is commonly used for porous electrode modeling at the macroscopic level \cite{bard2022electrochemical, bui2022continuum, wlodarczyk2023continuum}, where spatiotemporal variables are represented using ``mean-field" approximations. In this framework, the electrode and electrolyte potential domains are superimposed, resulting in a system of coupled Poisson equations governed by charge conservation and electrochemical kinetics. The corresponding coupled potential field model for the porous electrode system is typically written in the following closure form \cite{you2009simple, krishnamurthy2011computational, bayanov2011numerical}:  
\begin{subequations} 
\begin{align}
\nabla \cdot (-\sigma \nabla \phi_e) &= \nabla \cdot \vec{j_e} \quad \text{in } \Omega_e, \label{eq:poisson-1} \\
\nabla \cdot (-\kappa \nabla \phi_l) &= \nabla \cdot \vec{j_l} \quad \text{in } \Omega_l, \label{eq:poisson-2} 
\end{align}
\end{subequations}
\aftergroup\noindent
where the two Poisson equations are nonlinearly coupled through charge conservation and reaction kinetics \cite{dickinson2020butler, alotto2014redox}, given by:
\setlength{\jot}{10pt} 
\begin{subequations} 
\begin{align}
&\nabla \cdot \vec{j_e} = -\nabla \cdot \vec{j_l} = -sj \label{eq:charge}\,, \\
&j = j_0 \left[ \exp\left(\frac{(1-\alpha) F \eta}{RT}\right) - \exp\left(-\frac{\alpha F \eta}{RT}\right) \right] \label{eq:b-v}\,, \\
&\eta = \phi_e - \phi_l - E_{eq} \,.\label{eq:eta}
\end{align}
\end{subequations}

\mwang{The descriptions to the symbols} in the above equations are summarized in Table \ref{tab:symbols}. Notably, the right-hand sides of \eqref{eq:poisson-1} and \eqref{eq:poisson-2} have opposite signs, canceling each other when summed. Moreover, they exhibit a trigonometric hyperbolic dependence on \mwang{the pair} \((\phi_e, \, \phi_l)\), or more precisely, on their difference \(\phi_e - \phi_l\). \mwang{Furthermore}, the current density \(j\) is negative for the reduction process and positive for the oxidation process.

\begin{table}[h!]
\centering
\caption{Description of symbols in \eqref{eq:poisson-1}--\eqref{eq:poisson-2} and \eqref{eq:charge}--\eqref{eq:eta}}
\renewcommand{\arraystretch}{1.2} 
\resizebox{\linewidth}{!}{ 
\small 
\begin{tabular}{p{4.5cm} p{9cm} p{2cm}}
\Xhline{0.6pt} 
\textbf{Symbol} & \textbf{Description} & \textbf{Unit} \\ \Xhline{0.6pt} 
\multicolumn{3}{l}{\textit{Primary unknowns:}} \\
\(\phi_e\)      & Porous electrode potential                       & V \\
\(\phi_l\)      & Electrolyte potential                            & V \\ \hline
\multicolumn{3}{l}{\textit{Secondary unknowns:}} \\
\(\vec{j_e}\)   & Current density in the porous electrode          & A/m\(^2\) \\
\(\vec{j_l}\)   & Current density in the electrolyte               & A/m\(^2\) \\
\(j\)           & Reaction current density                         & A/m\(^2\) \\ 
\(\eta\)        & Overpotential                                    & V \\ \hline
\multicolumn{3}{l}{\textit{Parameters:}} \\
\(\sigma\)      & Electrical conductivity of the porous electrode  & S/m \\
\(\kappa\)      & Ionic conductivity of the electrolyte            & S/m \\ 
\(\alpha\)      & Electron transfer coefficient                    & - \\
\(s\)           & Specific surface area of the porous electrode    & m\(^{-1}\) \\
\(j_0\)         & Exchange current density                         & A/m\(^2\) \\
\(F\)           & Faraday constant, \(96485\)                     & \(\mathrm{C/mol}\) \\
\(R\)           & Ideal gas constant, \(8.314\)                    & \(\mathrm{J/(mol \cdot K)}\) \\
\(T\)           & Operating temperature                            & K \\ 
\(E_{eq}\)      & Equilibrium potential                            & V \\
\(\Omega_e\)    & The modeling domain of electrode                 & - \\
\(\Omega_l\)    & The modeling domain of electrolyte               & - \\
\Xhline{0.6pt} 
\end{tabular}
}
\label{tab:symbols}
\end{table}

Although it has not been framed this way in the relevant literature, \eqref{eq:poisson-1}--\eqref{eq:poisson-2} and \eqref{eq:charge}--\eqref{eq:eta} describe a dual-continuum approach in macroscopic modeling \cite{barenblatt1960basic, wang2024physical}. This approach assumes that the two coupled dynamics are spatially superimposed, with their interaction governed by a pair of antisymmetric source/sink terms that enforce global conservation. By abstracting the system in this way, the dual-continuum approach effectively alleviates the need to resolve local porous geometry, significantly improving computational efficiency. This modeling strategy, along with its generalized versions, has been extensively studied in the porous media flow community as part of multiscale frameworks designed to upscale the coexistence of fast and slow fluid flow dynamics in complex media with well-distributed high- and low-conductivity regions \cite{gerke1993dual, arbogast1997computational, akkutlu2017multiscale, wang2020generalized, vogel2010physical}.

In the context of porous electrodes, the dual-continuum setup considers the solid electrode and liquid electrolyte phases as coexisting everywhere and geometrically superimposed, such that \(\Omega_e = \Omega_l\). The porous electrode system exhibits two distinct conductivity regimes in the sense that the solid phase has significantly higher electrical conductivity than the electrolyte phase. The right-hand sides of \eqref{eq:poisson-1} and \eqref{eq:poisson-2} form an antisymmetric pair, which can also be interpreted as internal or implicit boundary (interfacial) conditions governing the interaction between electrode and electrolyte dynamics. This coupling is generally described by the phenomenological Butler-Volmer equation \cite{dickinson2020butler, alotto2014redox, guidelli2014defining}, along with some forms of correction for surface species concentrations \cite{dreyer2016new, seidenberg2025interpreting}. Unlike the porous media fluid flow applications, where such interfacial conditions are often approximated as pseudo-steady \cite{gerke1993dual, arbogast1997computational, akkutlu2017multiscale, wang2020generalized, vogel2010physical}, the electrochemical system exhibits strong nonlinearity, following a hyperbolic semi-\(\sinh\) dependence. 

The numerical solution of \eqref{eq:poisson-1} and \eqref{eq:poisson-2} presents challenges primarily in two aspects. First, the right-hand side is highly nonlinear, exhibiting a hyperbolic dependence on the unknowns (\(\phi_e\) and \(\phi_l\)). Even a small variation in (\(\phi_e - \phi_l\)) can lead to significant changes in the right-hand side, potentially causing numerical instability. Second, a major challenge arises from the boundary conditions, which are intricate due to the coupled nature of the system and are dictated by the operating mode of the electrochemical cell. Electrochemical cells are typically operated under either galvanostatic or potentiostatic control, with galvanostatic operation being more common \cite{bard2022electrochemical, newman2021electrochemical}. Each mode determines the boundary conditions used in modeling. In galvanostatic operation, the current flux through the cell is fixed and externally controlled, allowing the potential to adjust dynamically in response to internal electrochemical reactions and mass transport. This applied current flux leads to an all-Neumann boundary setup for the two Poisson equations, rendering the system under-constrained and making standard numerical solution strategies unsuitable due to singularity issues. In contrast, under potentiostatic operation, the potential difference between the electrodes (typically between the working and reference electrodes) is maintained at a prescribed value, while the current adapts accordingly to accommodate reaction kinetics and transport processes. This results in a Dirichlet boundary condition for the electrode potential, providing a much more well-posed setup for numerical solutions. 

Most continuum-level porous-electrode studies rely on commercial solvers to compute galvanostatic responses of the coupled system and report simulation results under homogeneous conductivity fields \cite{krishnamurthy2011computational, xu2015fundamental, gandomi2015situ, chen2023hybrid, shah2008dynamic, martinez2024computational, you2009simple}. These works typically present the governing equations and all-Neumann boundary conditions, but give little discussion of the numerical treatment needed to obtain a well-posed solution. As a result, the computational methodology is difficult to reproduce, assess, or transfer to heterogeneous settings and alternative numerical frameworks. This work is thus motivated to provide the explicit numerical schemes for the galvanostatic case and a basis for implementation and benchmarking beyond black-box solvers.
 
The main contributions of this paper are as follows:  
(1) we prove that, under galvanostatic conditions, the coupled system admit a unique solution in terms of \((\phi_e - \phi_l)\) and a unique pair of \((\phi_e ,\, \phi_l)\), up to a shared constant shift in the solution space of \(\phi_e\) and \(\phi_l\).  
(2) we provide two methods to eliminate this constant shift: a Lagrange-multiplier constraint derived from a unified energy functional, and a Gauss-theorem-based Dirichlet substitution.
(3) we develop two nonlinear solution schemes for solving the galvanostatic problem: a decoupled approach and a fully coupled approach.  
(4) we demonstrate that, without introducing any additional reference potential, a solution can still be obtained using global constraining in a purely mathematical sense to address the singular system in the fully coupled scheme.  

\section{Mathematical Model} \label{sec:model}
\vspace{-1.8ex}
\subsection{Governing Equations} \label{subsec:governing equations}
\vspace{-1.8ex}
We consider a symmetric single-electron redox reaction and set \(\alpha=0.5\) \cite{bard2022electrochemical,marcus1997electron}. We focus on reduction (\(j<0\)) \cite{newman2021electrochemical}. The methods in Sections~\ref{sec:methods-galvanostatic} and~\ref{sec:methods-potentiostatic} extend to oxidation, asymmetric \(\alpha\), and transient transport.
With these assumptions, \eqref{eq:poisson-1}--\eqref{eq:eta} reduce to
\begin{subequations}\label{eq:poisson_reduced}
\begin{align}
\nabla\!\cdot(-\sigma\nabla\phi_e) &= -a\,\sinh(b\eta) \quad \text{in } \Omega_e, \label{eq:poisson-11}\\
\nabla\!\cdot(-\kappa\nabla\phi_l) &= \phantom{-}a\,\sinh(b\eta) \quad \text{in } \Omega_l, \label{eq:poisson-22}
\end{align}
\end{subequations}
where \(\eta:=\phi_e-\phi_l-E_{eq}\), \(a=2sj_0\), and \(b=\tfrac{0.5F}{RT}\).
\subsection{Model Geometry}
\begin{figure} [ht!]
    \centering
    \includegraphics[width=0.4\linewidth]{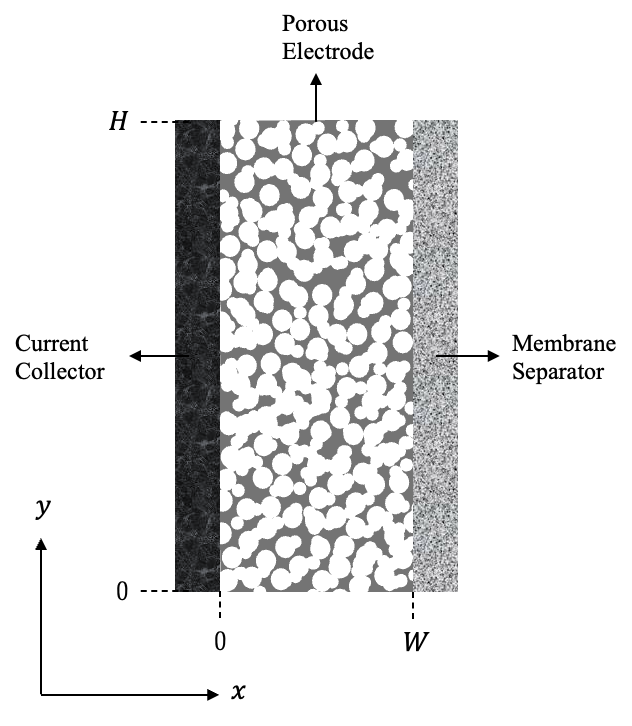}
    \caption{The half-cell porous electrode model geometry} 
    \label{fig:model-geometry}
\end{figure}

\vspace{-1.8ex}
The model geometry considered in this paper is illustrated in Figure~\ref{fig:model-geometry}. The domain representing a half-cell setup comprises a porous electrode positioned between a current collector on the left and a membrane separator on the right. The porous electrode is modeled as a composite material comprising a solid matrix (shown in gray) and pore spaces (shown in white), the latter being filled with electrolyte. 
The model extends over a width \(W\) in the \(x\)-direction and a height \(H\) in the \(y\)-direction, with the current collector and membrane separator defining the left and right boundaries, respectively. The entire system can be characterized by its dual-continuum nature, consisting of a solid electrode-continuum \(\Omega_e\) and a liquid-continuum \(\Omega_l\) that are spatialy superimposed.
In the following, when we use the symbol \(\Omega\), it denotes both the electrode and the electrolyte domains without explicitly distinguishing them. The redox reaction interface between these two continua is formed by the surface of the pore space. It acts as an implicit internal interfacial boundary where electrochemical reactions occur. This dual-continuum configuration serves as the foundation for describing the coupled potential fields arising from electronic and ionic conduction.

\subsection{Boundary Conditions}
\vspace{-1.8ex}
Here, we consider both the galvanostatic and potentiostatic operating conditions. We assume the boundaries consist entirely of either Dirichlet or Neumann conditions.

\textbf{Galvanostatic.}
The left boundary \(x = 0\) of \(\Omega_e\), attaching to the current collector, is subjected to a fixed current density applied through the external circuit. Since the membrane separator is not conductive to electrons, the right boundary \(x = W\) of \(\Omega_e\) is considered to have no electron current flux. Additionally, the top and bottom boundaries are assumed to be insulated to current flux. The boundary conditions for the porous electrode are:
\begin{equation}
\sigma \frac{\partial \phi_e}{\partial n} =
\begin{cases} 
    -j_{\text{applied}}, & \text{at } x = 0, \\
    \addlinespace[12pt]  
    0, & \text{at } x = W, \, y = 0, \, y = H,
\end{cases}
\label{eq:BCe-g}
\end{equation}
where \(j_{\text{applied}}\) is the applied current density (\(A/m^2\)) and \(\frac{\partial \phi_e}{\partial n}\) denotes the directional derivative in the outward normal to the boundary. 
For the electrolyte domain, its left boundary \(x = 0\) does not support ionic conduction, resulting in no ionic current flux. Conversely, at its right boundary \(x = W\) in contact with the membrane separator, there is supporting ionic flux between the positive and negative electrolyte regions. This ionic current flux balances the charge changes caused by electrochemical reactions and completes the internal current circuit. Under galvanostatic operation, the ionic current density at \(x=W\) matches the current density applied to the electrode via the current collector. The top and bottom boundaries are insulated with respect to ionic flux. We then have for the electrolyte:
\setlength{\jot}{10pt} 
\begin{equation}
\kappa \frac{\partial \phi_l}{\partial n} = 
\begin{cases} 
    j_{\text{applied}}, & \text{at } x = W, \\
    \addlinespace[12pt]  
    0, & \text{at } x = 0, \, y = 0, \, y = H.
\end{cases}
\label{eq:BCl-g}
\end{equation}
Here \(j_{\text{applied}}\) represents a non-negative value.

\textbf{Potentiostatic.}
The left boundary \(x = 0\) of \(\Omega_e\), attaching to the current collector, is set to an externally specified potential as a Dirichlet boundary condition. The remaining boundaries are the same as the galvanostatic case. Thus, we have:
\setlength{\jot}{10pt}
\begin{equation}
\begin{cases}
\phi_e  =  V_{\text{applied}},&\mwang{\text{at }x = 0,} \\
    \addlinespace[12pt]  
\sigma \frac{\partial \phi_e}{\partial n} = 0,&\text{at }x = W, \, y = 0, \, y = H. 
\end{cases}
\label{eq:BCe_p}
\end{equation} 

For the electrolyte domain, the boundaries at \(x = 0\), \(y = 0\), and \(y = H\) of \(\Omega_l\) remain under zero-flux conditions. However, the boundary at \(x = W\) is not  explicitly specified by the operating condition. We can impose either Dirichlet or Neumann boundaries there and then numerically sweep the respective boundary values when solving the system. Here we consider two cases:

\setlength{\jot}{10pt}
\begin{equation}
\begin{cases}
\phi_l = V_{\text{sweep}}, & \text{at }x = W,\\
    \addlinespace[12pt]  
\kappa \frac{\partial \phi_l}{\partial n} = 0, &\text{at }x = 0, \, y = 0, \, y = H,
\end{cases}
\label{eq:BCl_p_dirichlet}
\end{equation}

\setlength{\jot}{10pt}
\begin{equation}
\kappa \frac{\partial \phi_l}{\partial n} = 
\begin{cases} 
    j_{\text{sweep}}, & \text{at } x = W, \\
    0, & \text{at } x = 0, \, y = 0, \, y = H,
\end{cases}
\label{eq:BCl_p_neumann}
\end{equation}

where \(V_{\text{sweep}}\) and \(j_{\text{sweep}}\) are to be swept during solution. 

\section{Numerical Methods: Galvanostatic Solution} \label{sec:methods-galvanostatic}
\vspace{-1.8ex}
Due to the nonlinear source/sink terms in the coupled Poisson system, we need to devise proper iterative linearization schemes for the galvanostatic solution.
However, the major challenge arises from the all-Neumann boundary conditions, which leads to non-uniqueness in the respective solution spaces of \(\phi_e\) and \(\phi_l\). Such a situation results in a singular system matrix after discretization. In this section, we provide decoupled and coupled numerical schemes to address this issue.
\vspace{-1.8ex}
\subsection{Solution Existence and Uniqueness} \label{sec:existence}
\vspace{-1.8ex}
\textbf{Existence.} 
\label{subsec:existence}
The solution existence requires satisfying a compatibility condition based on Gauss's theorem, as shown in \eqref{eq:Gauss_e} and \eqref{eq:Gauss_l}. In a physical sense, this requirement is interpreted as charge conservation. 
\setlength{\jot}{10pt}
\begin{subequations} 
\begin{align}
    \int_{\Omega_e} (-f) \, d\Omega_e &= -\int_{\partial \Omega_e} \sigma \frac{\partial \phi_e}{\partial n} \, dS_e, \label{eq:Gauss_e} \\
    \int_{\Omega_l} f \, d\Omega_l &= -\int_{\partial \Omega_l} \kappa \frac{\partial \phi_l}{\partial n} \, dS_l. \label{eq:Gauss_l}
\end{align}
\end{subequations}

Here, \( f = a \, \sinh(b \, \eta) \), while \( \Omega_{e/l} \) and \( \partial \Omega_{e/l} \) represent the interior and boundary of the model domain \( \overline{\Omega}_{e/l} \), respectively. \( dS_{e/l} \) denotes a differential surface element on \( \partial \Omega_{e/l} \).

\textbf{Uniqueness.}
For this nonlinearly coupled all-Neumann system, the potentials are not individually unique. The solution set consists of pairs \((\phi_e,\phi_l)\) that are unique up to a shared additive constant (see \textbf{SI.A}),
\begin{equation}
\{(\phi_e,\phi_l)\mid \phi_e=\phi_e^*+C,\ \phi_l=\phi_l^*+C,\ C\in\mathbb{R}\},
\label{eq:unique_pair}
\end{equation}
where \((\phi_e^*,\phi_l^*)\) is any particular solution pair satisfying \eqref{eq:poisson-11}--\eqref{eq:poisson-22} with boundary conditions \eqref{eq:BCe-g}--\eqref{eq:BCl-g}. \eqref{eq:unique_pair} immediately yields the unique potential difference,
\begin{equation}
\phi_e-\phi_l=\phi_e^*-\phi_l^*,
\label{eq:unique_pair_}
\end{equation}
and explains why the galvanostatic case is numerically delicate: \(\phi_e\) and \(\phi_l\) can drift by a constant unless one additional reference value is imposed to remove the nullspace. In practice, this can be achieved by fixing \(\phi_e(x^*)\) or \(\phi_l(x^*)\) at an arbitrary location \(x^*\). For instance, we ground a point in the electrode or the current collector boundary/near-boundary region). Then the remaining task is to enforce this reference consistently within the chosen nonlinear solution strategy; below we present two ways to impose such a fixed-location reference and give decoupled and fully coupled numerical schemes.

\subsection{Constant Potential Referencing} \label{sec:constant-referencing}
\vspace{-1.8ex}
\textbf{Lagrange Constraint Method (LCM).} \label{sec:LCM}
To fix the constant-shift degree of freedom, we impose a fixed-location potential reference in \(\Omega_e\) via Lagrange multipliers, while keeping the galvanostatic boundary conditions unchanged. Although the reference location and value are arbitrary, we set \(\phi_e=0\) at selected points near the left boundary of \(\Omega_e\) as a physically motivated choice. For the coupled system, a unified energy functional yields (see \textbf{SI.B}) the augmented formulation with boundary conditions in \eqref{eq:BCe-g}--(\eqref{eq:BCl-g}:
\begin{subequations}
\begin{align}
\nabla \cdot (-\sigma \nabla \phi_e) + f + \sum_{i=1}^{m}\lambda_i\,\delta(x-x_i) &= 0 \quad \text{in } \Omega_e, \label{eq:extended-lagrange-1_}\\
\nabla \cdot (-\kappa \nabla \phi_l) - f &= 0 \quad \text{in } \Omega_l, \label{eq:extended-lagrange-2_}\\
\phi_e(x_i) - c_i &= 0 \quad \text{for } i=1,\dots,m, \label{eq:extended-lagrange-3_}
\end{align}
\end{subequations}
where $x_i\in\Omega_e$ specify the constraint locations $\{x_1,\dots,x_m\}$, $\delta(\cdot)$ denotes the Dirac delta, and $\lambda_i$ are the associated Lagrange multipliers; we collect them as $\Lambda_e:=\{\lambda_i\}_{i=1}^m$. Thus, \(\phi_e(x_i)\) is pinned to \(c_i\) and serves as the reference for the coupled galvanostatic solution. An analogous construction applies if the reference is imposed on \(\phi_l\). Other local constraint techniques include penalty terms \cite{zhu1998modified} and direct matrix modification \cite{quarteroni1994numerical}: the former enforces the constraint only approximately and can be sensitive to its heuristic penalty parameter, while the latter enforces $\phi=C$ by replacing one row of the discrete system, thereby discarding the original governing equation at that location and potentially creating a local imbalance in the discrete flux/residual; In contrast, the Lagrange approach enforces the constraint exactly in a theoretically consistent manner, preserving the original PDE structure (see \textbf{SI.B}), but requiring a larger algebraic system.

\textbf{Dirichlet Substitution Method (DSM).} \label{sec:DSM}
Recall Section~\ref{subsec:existence}. By Gauss's theorem, a Poisson equation satisfies a global compatibility condition. The net boundary flux over $\partial\Omega$ equals the total source/sink over $\Omega$. This holds for both Neumann and Dirichlet problems, since Dirichlet boundaries imply a flux even when it is not prescribed. Therefore, in the galvanostatic setup, we may replace one Neumann boundary by an arbitrary Dirichlet reference potential without changing the original coupled system, because the substituted boundary flux is still determined by the same global balance together with the remaining Neumann boundaries. For example, under the galvanostatic boundary conditions in \eqref{eq:BCe-g}--\eqref{eq:BCl-g}, we may replace the electrode Neumann condition at $x=0$ by an arbitrary Dirichlet reference $\phi_e|_{x=0}=\phi_{\mathrm{ref}}$. The flux across the substituted boundary is not an additional free input; it is determined implicitly by the same global balance and therefore matches the flux that would arise in the original all-Neumann formulation. In particular,
\begin{equation}
\int_{\partial \Omega_{e,\mathcal{D}}\,|_{x=0}} \sigma \frac{\partial \phi_e}{\partial n}\, dS_e
= -\int_{\Omega_e} (-f)\, d\Omega_e
= \int_{\Omega_l} f\, d\Omega_l
= -\int_{\partial \Omega_{l,\mathcal{D}}\,|_{x=W}} \kappa \frac{\partial \phi_l}{\partial n}\, dS_l,
\label{eq:dirichlet-replace}
\end{equation}
where $\mathcal{N}$ and $\mathcal{D}$ denote the Neumann and Dirichlet portions of the boundary, respectively.

\subsection{Nonlinear Solution Schemes}
\vspace{-1.8ex}
In addition to imposing a constant reference potential, it is essential to establish internal referencing between \(\phi_e\) and \(\phi_l\) during the nonlinear solution process.
We consider both decoupled schemes, where the two equations are solved sequentially, and fully coupled schemes, where they are solved simultaneously.
\vspace{-1.8ex}
\subsubsection{Decoupled Scheme} \label{sec:fixed-pint}
\vspace{-1.8ex}
\mwang{We use a staggered decoupling scheme for the coupled system \eqref{eq:poisson-11}--\eqref{eq:poisson-22}. At iteration $k+1$, we update $\phi_e$ and $\phi_l$ sequentially. Each solve uses the latest available information from the other domain at iteration $k$. This turns the coupled problem into a sequence of independent subproblems. Table~\ref{tab:picard-steps} summarizes the procedure. In \eqref{eq:poisson-111} and \eqref{eq:poisson-222}, choosing $m=k$ gives a fully explicit scheme, whereas $m=k+1$ requires an implicit solve. We use $m=k+1$ because the explicit option is generally unstable for elliptic problems.} \mwang{Our stopping criterion is designed for the strong sensitivity of $\sinh(\eta)$ to changes in $\eta$, which provides a reliable convergence monitor.} This decoupling raises two issues: the coupling between $\phi_e$ and $\phi_l$ is only weakly imposed through the source term $f$, since the sequential solves do not enforce a direct gridblock-wise referencing between the two potentials; and imposing a constant reference in $\Omega_e$ removes the singularity in \eqref{eq:poisson-111}, but \eqref{eq:poisson-222} remains singular. To address both, we introduce an additional reference potential $c_l$ in $\Omega_l$. This resolves the singular matrix in \eqref{eq:poisson-222} and provides a cross-referencing mechanism within the sequential scheme. However, \eqref{eq:unique_pair} shows that only one global reference is needed to isolate a unique solution pair, so $c_l$ cannot be prescribed after $\Omega_e$ is anchored. We therefore treat $c_l$ as an unknown and determine it from a closure condition. We use the compatibility constraint, or charge conservation, in \eqref{eq:Gauss_e}--\eqref{eq:Gauss_l}). In practice, $c_l$ can be found using a search algorithm. Below, we show that both LCM and DSM can be adapted to iteratively solve for $c_l$ in a self-consistent manner.
\begin{table}[h!]
\centering
\caption{Decoupled nonlinear solution scheme}
\renewcommand{\arraystretch}{1.3} 
\resizebox{\linewidth}{!}{ 
\small 
\begin{tabular}{p{3cm} >{\raggedright\arraybackslash}p{12cm}}
\Xhline{0.6pt} 
\textbf{Step} & \textbf{Description} \\ 
\Xhline{0.6pt} 
\textit{Initialization} & 
Set the initial guess for the nonlinear term \(f^0(\eta)\) using average local current density. \\ \hline
\textit{Iteration} & 
\vspace{-1.3\baselineskip} 
\begin{enumerate}[left=0pt]
    \item Solve the following two Poisson equations sequentially: 
        \begin{subequations} 
        \setlength{\jot}{10pt} 
        \begin{gather}
            \nabla \cdot (-\sigma \, \nabla \phi_e^{k+1})  = -f(\phi_e^{m} - \phi_l^{k} - E_{\text{eq}}), \label{eq:poisson-111} \\
            \nabla \cdot (-\kappa \, \nabla \phi_l^{k+1}) = f(\phi_e^{k} - \phi_l^{m} - E_{\text{eq}}). \label{eq:poisson-222} 
        \end{gather}
        \end{subequations}
    subject to boundary conditions \eqref{eq:BCe-g} and \eqref{eq:BCl-g}. \(k\) is the iteration index. 
    \item Update the nonlinear term to \mwang{\(f^{k+1}:= f(\eta^{k+1})\)} for convergence check letting 
    \(\eta^{k+1} = \phi_e^{k+1} - \phi_l^{k+1} - E_{\text{eq}}.\)
    \item Check convergence: If \(\|f^{k+1} - f^k\|\) falls below a predefined threshold, proceed to output. Otherwise, continue the update.
\end{enumerate} \\ \hline
\textit{Output} & 
Return the converged \(\phi_e\) and \(\phi_l\). \\ 
\Xhline{0.6pt} 
\end{tabular}
}
\label{tab:picard-steps}
\end{table}

\textbf{LCM Discretization.}
We first write the discretized form of \eqref{eq:extended-lagrange-1_} -- \eqref{eq:extended-lagrange-3_}, which reads 
\setlength{\jot}{10pt} 
\begin{subequations}
\begin{align}
    A_e \Phi_e + \Delta_e^T \Lambda_e &= F_e, \label{eq:poisson_lagrange_discretization} \\
    \Delta_e \Phi_e &= C_e, \label{eq:constraint_discretization} 
\end{align}
\end{subequations}
where \(A_e\) is the system matrix, \(F_e\) is the discretized right-hand side vector, \(\Phi_e := [\phi_1^{k+1}, \phi_2^{k+1}, \dots]^T\), \(\Lambda_e\) is the Lagrange multiplier vector, and \(C_e: = c_e\,\mathbf{1}_{n_c}\) is the reference potential vector, where $\mathbf{1}_{n_c} =[1,1,\cdots,1]^T\in\mathbb{R}^{n_c}$ and  \(n_c\) is the number of constrained gridblocks. The Neumann boundary condition information is incorporated \(A_e\) and \(F_e\). \(\Delta_e\) is a Dirac delta matrix with entries of \(1\) at the locations corresponding to the constraints, and \(0\) elsewhere. It projects the \(\Lambda_e\) and \(\Phi_e\) vectors onto the indices of the constrained gridblocks. The matrix \(\Delta_e\) has dimensions \(n_c \times N\), where \(N\) is the total number of gridblocks. Each row of this sparse matrix corresponds to a constrained gridblocks. The corresponding discretization of \(\Omega_l\) can be done in an analogous way. The resulting linear systems are:
\begin{subequations}
\begin{align}
\begin{bmatrix}
    A_e & \Delta_e^T \\
    \Delta_e & 0
\end{bmatrix}
\begin{bmatrix}
    \Phi_e \\
    \Lambda_e
\end{bmatrix}^{k+1}
&=
\begin{bmatrix}
    F_e^{k/k+1} \\
    C_e
\end{bmatrix}, \label{eq:lagrange_matrix_e} \\[10pt]
\begin{bmatrix}
    A_l & \Delta_l^T \\
    \Delta_l & 0
\end{bmatrix}
\begin{bmatrix}
    \Phi_l \\
    \Lambda_l
\end{bmatrix}^{k+1}
&=
\begin{bmatrix}
    F_l^{k/k+1} \\
    C_l^{*}
\end{bmatrix}. \label{eq:lagrange_matrix_l}
\end{align}
\end{subequations}
where \(C_l^*: = c_l^*\,\mathbf{1}_{n_l}\) is an \mwang{vector parametrized by an unknown $c_l^{*}$}. The superscript \(k/k+1\) indicates, for example, \(F_e\) is evaluated using \mwang{\(\phi_e^{k/k+1}\).}

\vspace{-1.8ex}
\textbf{DSM Discretization.}
In this case, the discretization is simpler as \(c_e\) and \(c_l^*\) are imposed as constant Dirichlet boundaries at \(x=0\) and \(x=W\). The resulting discretized system is given by:
\begin{equation}
A_e\, \Phi_e^{k+1} = F_e^{k/k+1} \quad \text{and} \quad A_l\, \Phi_l^{k+1} = F_l^{k/k+1}
\end{equation}
where the corresponding Dirichlet boundaries are incorporated in (\(A_e\), \(F_e\)) and  (\(A_l\), \(F_l\)). 

\textbf{Unknown Reference Potential Search.} \label{sec:search-algorithm}
We determine the unknown reference potential \(c_l^*\) by formulating an optimization problem to isolate its unique value that satisfies \eqref{eq:Gauss_e}and \eqref{eq:Gauss_l}. Together with \eqref{eq:dirichlet-replace}, we define the following objective function: 
\begin{equation}
    c_l^* = \arg \min_{c_l} \left\| \int_{\mathcal{B}} j_{\text{applied}} \, d{S} +\int_\Omega f(\cdot;c_l) \, d\Omega \right\|^2,\quad
    \mathcal{B} = \partial \Omega|_{x=0} \,\, \text{or} \,\,\mathcal{B}=\partial \Omega|_{x=W}
    . \label{eq:objective_function}
\end{equation}
To solve this optimization problem, we employ a gradient descent approach with a line search mechanism \cite{bonnans2006numerical} that dynamically adjusts the search step size (see \textbf{SI.C}). This method iteratively refines the reference potential with linear convergence.
In Table \ref{tab:line-search-algorithm}, we outline the complete solution strategy. It incorporates the staggered iteration as an inner loop for solving the nonlinear Poisson equations sequentially, with either the LCM or DSM approach used to impose \(c_e\) and \(c_l^*\). The combination of the outer search loop and the inner Picard iteration thus provides a way for obtaining galvanostatic solutions in a decoupled fashion while ensuring internal referencing between \(\phi_e\) and \(\phi_l\). This approach can address the singularity issues and enforces charge conservation in a self-consistent manner.

\begin{table}[H]
\centering
\caption{Galvanostatic solution strategy using decoupled scheme}
\renewcommand{\arraystretch}{1.3} 
\resizebox{\linewidth}{!}{ 
\small 
\begin{tabular}{p{3cm} >{\raggedright\arraybackslash}p{12cm}}
\Xhline{0.6pt} 
\textbf{Step} & \textbf{Description} \\ 
\Xhline{0.6pt} 
\textit{Initialization} & 
Set a initial reference potential \(c_l^*\).
 \\ \hline
\textit{Solve} & 
Refer to Table \ref{tab:picard-steps}. Use LCM or DSM to solve the coupled nonlinear Poisson equations in a sequential fashion. Iterate until convergence or reaching stopping criteria. \\ \hline
\textit{Optimization} & 
\vspace{-1.3\baselineskip} 
\begin{enumerate}[left=0pt]
    \item Check convergence using \eqref{eq:objective_function}: If \mwang{the mismatch is below the prescribed tolerance}, exit and output the solution \(\phi_e\) and \(\phi_l\); otherwise, continue.
    \item Compute the gradient of \eqref{eq:objective_function}.
    \item Perform line search to obtain an optimized step size.
    \item Update \(c_l^*\) and \mwang{rerun to solve the decoupled PDEs.}
\end{enumerate} \\ \hline
\textit{Output} & 
Return the globally converged solution \(\phi_e\) and \(\phi_l\). \\ 
\Xhline{0.6pt} 
\end{tabular}
}
\label{tab:line-search-algorithm}
\end{table}

\subsubsection{Fully Coupled Scheme} \label{newton-raphson}
\vspace{-1.8ex}
A fully coupled formulation can enforce gridblock-wise cross-referencing between $\phi_e$ and $\phi_l$, which removes the need to introduce and search for $c_l^*$. We can solve \eqref{eq:poisson-11}--\eqref{eq:poisson-22} simultaneously using a fully coupled Jacobian system, treating $\phi_e$ and $\phi_l$ implicitly in the nonlinear solve. In contrast to the decoupled scheme, which holds one field fixed while updating the other, the coupled Jacobian captures the cross-dependence of the two potentials and suppresses arbitrary constant shifts by enforcing internal referencing within the assembled system. We demonstrate both LCM and DSM to impose a constant reference potential, and we also present a global constraining method that first computes the unique difference $\phi_e-\phi_l$ and then recovers $\phi_e$ and $\phi_l$ by post-processing under an arbitrary constant reference.

\textbf{LCM Discretization.} 
Discretizing \textbf{Eqs.~(\ref{eq:extended-lagrange-1_})--(\ref{eq:extended-lagrange-3_})} gives:
\begin{subequations}
\begin{align}
    Q_e^L &= A_e \Phi_e - F_e + \Delta_e^T \Lambda_e, \label{eq:residual_e_lagrange} \\
    Q_l^L &= A_l \Phi_l - F_l, \label{eq:residual_l_lagrange} \\
    G_e^L &= \Delta_e \Phi_e - C_e. \label{eq:residual_e_larange_constraint} 
\end{align}
\end{subequations}
where $F_e$ and $F_l$ are the discrete forms of $-f$ and $f$, and the superscript $L$ denotes Lagrange. Define the stacked residual and unknown as $R^L=[Q_e^L,\,Q_l^L,\,G_e^L]^T$ and $U^L=[\Phi_e,\,\Phi_l,\,\Lambda_e]^T$.
The extended Jacobian system is
\begin{equation}
J^L \tau_U^L = -R^L,
\qquad
J^L = \frac{\partial R^L}{\partial U^L},
\label{eq:jacobian_lagrange_system}
\end{equation}
with update $U^{L,k+1}=U^{L,k}+\tau_U^L$. The Jacobian blocks are
\begin{equation}
\begin{aligned}
\frac{\partial Q_e^L}{\partial \Phi_e} &= A_e - \frac{\partial F_e}{\partial \Phi_e}, \qquad
\frac{\partial Q_e^L}{\partial \Phi_l} = - \frac{\partial F_e}{\partial \Phi_l}, \qquad
\frac{\partial Q_e^L}{\partial \Lambda_e} = \Delta_e^T,\\
\frac{\partial Q_l^L}{\partial \Phi_e} &= - \frac{\partial F_l}{\partial \Phi_e}, \qquad
\frac{\partial Q_l^L}{\partial \Phi_l} = A_l - \frac{\partial F_l}{\partial \Phi_l}, \qquad
\frac{\partial G_e^L}{\partial \Phi_e} = \Delta_e,
\end{aligned}
\label{eq:jacobian_lagrange_}
\end{equation}
so that
\begin{equation}
\renewcommand{\arraystretch}{2.2}
\setlength{\arraycolsep}{10pt}
\begin{bmatrix}
A_e - \frac{\partial F_e}{\partial \Phi_e} & -\frac{\partial F_e}{\partial \Phi_l} & \Delta_e^T\\
-\frac{\partial F_l}{\partial \Phi_e} & A_l - \frac{\partial F_l}{\partial \Phi_l} & 0\\
\Delta_e & 0 & 0
\end{bmatrix}
\begin{bmatrix}
\tau_{\Phi_e}\\ \tau_{\Phi_l}\\ \tau_{\Lambda_e}
\end{bmatrix}
=
-
\begin{bmatrix}
Q_e^L\\ Q_l^L\\ G_e^L
\end{bmatrix}.
\label{eq:jacobian_lagrange_matrix_system}
\end{equation}
Here $G_e^L$ includes the imposed reference potential. The off-diagonal blocks $-\partial F_e/\partial \Phi_l$ and $-\partial F_l/\partial \Phi_e$ encode the cross-coupling between $\phi_e$ and $\phi_l$, and the constraint rows remove the constant-shift nullspace, yielding a unique solution pair. It should be noted that, because the nonzero entries of $\Delta_e$ are $1$, the extended system can be poorly scaled. We recommend scaling the Lagrange components separately to match the magnitude of the PDE blocks before assembling the final Jacobian.

\textbf{DSM Discretization.} \label{sec:newton-raphson-dsm}
The Jacobian system of DSM can be readily obtained by removing the Lagrange components and applying appropriate modifications to \(A_e\) and \(A_l\), as well as \([Q_e]\) and \([Q_l]\), to accommodate the Dirichlet boundary conditions. The resulting Jacobian system is:
\begin{equation}
\renewcommand{\arraystretch}{3} 
\setlength{\arraycolsep}{10pt}    
\text{\normalsize$
\begin{bmatrix}
    A^d_e - \frac{\partial F_e}{\partial \Phi_e} & - \frac{\partial F_e}{\partial \Phi_l} \\
    - \frac{\partial F_l}{\partial \Phi_e} & A_l - \frac{\partial F_l}{\partial \Phi_l}
\end{bmatrix}
$}
\text{\normalsize$
\begin{bmatrix}
    \tau_{\Phi_e} \\ \tau_{\Phi_l}
\end{bmatrix}
$}
=
\text{\normalsize$
\begin{bmatrix}
    Q^{d,r}_e \\ Q_l 
\end{bmatrix}.
$}
\label{eq:jacobian_matrix_system_dirichlet}
\end{equation}
where the subscript \(d\) denotes the incorporation of Dirichlet boundary substitutions, and the superscript \(r\) indicates the inclusion of the constant reference potential.

\textbf{Global Constraining Method (GCM).} 
By \eqref{eq:unique_pair}--\eqref{eq:unique_pair_}, the difference $\phi_e-\phi_l$ is unique. If $\phi_e$ and $\phi_l$ are cross-referenced gridblock-wise, the original Jacobian system,
\begin{equation}
\renewcommand{\arraystretch}{3}
\setlength{\arraycolsep}{10pt}
\text{\normalsize$
\begin{bmatrix}
A_e - \frac{\partial F_e}{\partial \Phi_e} & - \frac{\partial F_e}{\partial \Phi_l} \\
- \frac{\partial F_l}{\partial \Phi_e} & A_l - \frac{\partial F_l}{\partial \Phi_l}
\end{bmatrix}
$}
\text{\normalsize$
\begin{bmatrix}
\tau_{\Phi_e} \\ \tau_{\Phi_l}
\end{bmatrix}
$}
=
\text{\normalsize$
\begin{bmatrix}
Q_e \\ Q_l
\end{bmatrix},
$}
\label{eq:jacobian_matrix_system_noref}
\end{equation}
can be used to determine a unique update for $\phi_e-\phi_l$ even though the system has a shared constant-shift nullspace. We then can target to resolve the singularity by imposing global constraining measures that select a stable solution among the infinitely many. We consider three measures: (i) $\min_{\tau_{\Phi}}\|\tau_{\Phi}\|^2$, which selects the minimum-norm update and removes the free constant shift \cite{ben2003generalized}; (ii) $\min_{\tau_{\Phi}}\|J\tau_{\Phi}-Q\|^2$, which seeks the least-squares residual solution and is commonly computed in a Krylov subspace \cite{golub2013matrix}; and (iii) $\min_{\tau_{\Phi}}\|J\tau_{\Phi}-Q\|^2+\lambda_T\|I\tau_{\Phi}\|^2$, which applies Tikhonov regularization \cite{golub1999tikhonov} and uses $\lambda_T$ to balance residual reduction and damping of $\tau_{\Phi}$. Here $\tau_{\Phi}=[\tau_{\Phi_e},\,\tau_{\Phi_l}]^T$ and $Q=[Q_e,\,Q_l]^T$. In practice, the Moore--Penrose pseudoinverse minimizes $\|\tau_{\Phi}\|^2$ \cite{barata2012moore}, MINRES minimizes $\|J\tau_{\Phi}-Q\|^2$ \cite{paige1975solution}, and a least-squares solver with Tikhonov regularization such as LSMR/LSQR minimizes $\|J\tau_{\Phi}-Q\|^2+\lambda_T\|I\tau_{\Phi}\|^2$ \cite{fong2011lsmr, paige1982lsqr}; we refer to the last option as LSTR (Least-Squares with Tikhonov Regularization). After obtaining any particular pair $(\phi_e^*,\phi_l^*)$, we recover the solution pair under an arbitrary constant reference $C$ by a simple shift. If we reference at $x_0$ with $\phi_e(x_0)=C$, the shift is $\phi_e^*(x_0)-C$ and the final pair is
\begin{equation}
(\phi_e,\phi_l)=(\phi_e^*,\phi_l^*)-[\phi_e^*(x_0)-C].
\label{eq:postprocessing}
\end{equation}

\section{Numerical Methods: Potentiostatic Solution} \label{sec:methods-potentiostatic}
\vspace{-1.8ex}
For the sake of completeness, we briefly outline the numerical strategy for solving the potentiostatic system. In comparison with the galvanostatic condition, the potentiostatic solution is straightforward, as the coupled system includes a fixed Dirichlet boundary condition for \(\Omega_e\) at \(x = 0\). Unlike the galvanostatic case, where boundary conditions of \(\Omega_l\) are explicitly specified at \(x = 0\) and \(x = W\), the boundary condition of \(\Omega_l\) at \(x = W\) in the potentiostatic system is not directly prescribed. Therefore, we consider two possible cases, as defined by \eqref{eq:BCl_p_dirichlet} and \eqref{eq:BCl_p_neumann}: Dirichlet at \(x = W\) and Neumann at \(x = W\). However, the system admits multiple solutions depending on the specific assigned boundary values at \(x = W\) for \(\Omega_l\), whether Dirichlet or Neumann. We thus can numerically sweep such boundary values to compute a series of different potentiostatic responses for both cases.

\textbf{Dirichlet at \(x = W\).} 
This case is trivial as the coupled system is well-constrained without any solution ambiguity. We can solve the system using either decoupled or fully-coupled approaches by sweeping the value \(V_{\text{sweep}} = \phi_l(x = W)\), provided that the sweeping boundary value remains within a physically meaningful and compatible range. For instance, in a charging process at the negative electrode, \(V_{\text{sweep}}\) should be larger than \(V_{\text{applied}}=\phi_e(x = 0)\).

\textbf{Neumann at \(x = W\).} 
The \(\Omega_l\) setup here is analogous to the galvanostatic situation. Since \(\phi_e(x = 0)\) is given, we can directly use \eqref{eq:jacobian_matrix_system_dirichlet} to solve this boundary flux sweeping problem. For the decoupled Picard iteration approach, we may need to introduce an unknown reference potential to \(\Omega_l\) using either the LCM or DSM method, and then iteratively determine this unknown along the way to obtain a self-consistent system solution pair \((\phi_e, \phi_l)\).

\section{Examples and Analysis} \label{sec:examples}
\vspace{-1.8ex}
To validate and analyze the introduced numerical methods, we consider examples based on the model shown in Figure~\ref{fig:model-geometry}. 
The examples model the charging process at a negative porous electrode, where the reduction reaction $\text{Ox}^{3+} + e^- \longrightarrow \text{Red}^{2+}$ converts a trivalent species $\text{Ox}^{3+}$ to a divalent species $\text{Red}^{2+}$ via a single-electron transfer. Table~\ref{tab:example-parameters} summarizes the physical and operational parameters. Here, $\sigma$ is corrected using a constant porosity of $0.78$ based on a solid electrode material conductivity of $1000~\mathrm{S/m}$. The electrolyte conductivity $\kappa$ is computed from the concentrations of $\text{Ox}^{3+}$ and $\text{Red}^{2+}$, set to $1053~\mathrm{mol/m^3}$ and $27~\mathrm{mol/m^3}$, respectively, and is corrected using $0.78$ porosity as well. The relevant correlations are available in \cite{krishnamurthy2011computational,bayanov2011numerical,shah2008dynamic}.
\begin{table}[ht!]
\centering
\caption{Physical and operational parameters used in the examples}
\renewcommand{\arraystretch}{1.4} 
\resizebox{\linewidth}{!}{ 
\small 
\begin{tabular}{p{4.5cm} p{10cm} p{3cm}} 
\Xhline{0.6pt} 
\textbf{Symbol} & \textbf{Description} & \textbf{Value} \\ 
\Xhline{0.6pt} 
\multicolumn{3}{l}{\textit{Electrode Geometry:}} \\
\(W\)           & Porous electrode width (m)                        & \(5 \times 10^{-3}\) \\
\(H\)           & Porous electrode height (m)                       & \(1 \times 10^{-1}\) \\
\(L\)           & Porous electrode length (m)                       & \(1 \times 10^{-1}\) \\ \hline
\multicolumn{3}{l}{\textit{Physical Parameters:}} \\
\(\sigma\)      & Electrical conductivity of the porous electrode (S/m)        & \(103.1891\) \\
\(\kappa\)      & Ionic conductivity of the electrolyte (S/m)                  & \(5.9514\) \\ 
\(s\)           & Specific surface area of the porous electrode (m\(^{-1}\))   & \(1.64 \times 10^4\) \\
\(\alpha\)      & Electron transfer coefficient (-)                    & \(0.5\) \\
\(j_0\)         & Exchange current density (A/m\(^2\))                 & \(2.7657\) \\
\(E_{eq}\)      & Equilibrium potential (V)                            & \(-0.1609\) \\
\(F\)           & Faraday constant (\(\mathrm{C/mol}\))                & \(96485\) \\
\(R\)           & Ideal gas constant (\(\mathrm{J/(mol \cdot K)}\))    & \(8.314\)  \\
\(T\)           & Operating temperature (K)                            & \(298.15\) \\ \hline
\multicolumn{3}{l}{\textit{Operational Parameters:}} \\
\(I_{\text{applied}}\)   & Applied galvanostatic current (A)    & \((0, \, 10]\) \\ 
\(V_{\text{applied}}\)   & Applied potentiostatic electrode potential (V) & \(0\) \\ 
\(V_{\text{sweep}}\)   & Sweeping potentiostatic electrolyte potential at \(x=W\) (V) & \([0.1, \, 0.5]\) \\ 
\(I_{\text{sweep}}\)     & Sweeping potentiostatic current at \(x=W\) (A) & \((0, \, 10]\) \\ 
\Xhline{0.6pt} 
\end{tabular}
}
\label{tab:example-parameters}
\end{table}

We focus on galvanostatic solutions, because the potentiostatic case uses standard numerical procedures. The main exception is when the boundary current is swept. In that setting, the problem can be viewed as a DSM-type variant of the galvanostatic solution. Aside from Section~\ref{sec:comparison_with_exact}, the remaining examples and analyses concern galvanostatic solutions. All numerical solutions, except for the global constraining cases, use BICGSTAB with an ILU preconditioner. 
\vspace{-1.8ex}
\subsection{Comparison with ``Exact" Solutions} \label{sec:comparison_with_exact}
\vspace{-1.8ex}
We compare the numerical solutions with their ``exact'' counterparts. We consider uniform $\sigma$ and $\kappa$ for analytical solutions in a one-dimensional setting. Thus, \eqref{eq:poisson-11} and \eqref{eq:poisson-22} become:
\begin{subequations}
\begin{gather}
\frac{d^2 \phi_e}{d x^2} = \frac{a}{\sigma} \sinh(b \eta), \label{eq:poisson-simplified1}\\
\frac{d^2 \phi_l}{d x^2} = -\frac{a}{\kappa} \sinh(b \eta), \label{eq:poisson-simplified2}
\end{gather}
\end{subequations}
where $a=2 s j_0$ and $b=\frac{0.5 F}{R T}$. Subtracting \eqref{eq:poisson-simplified2} from \eqref{eq:poisson-simplified1} gives
\(
\frac{d^2 \phi_e}{d x^2} - \frac{d^2 \phi_l}{d x^2}
= \frac{d^2(\phi_e-\phi_l)}{d x^2}
= \frac{d^2(\eta+E_{eq})}{d x^2}
= \frac{d^2 \eta}{d x^2},
\)
and thus \(\frac{d^2 \eta}{d x^2} = c \sinh(b \eta)\),
where $c=\left(\frac{a}{\sigma}+\frac{a}{\kappa}\right)$. Using the First Integral Method \cite{ghosh2021first}, we obtain the ``exact'' solution in implicit form (see \textbf{SI.D})
\begin{equation}
x = \int_{\eta_1}^{\eta}
\frac{d\eta}{\sqrt{\frac{2c}{b}\cosh(b\eta) + q_1^2 - \frac{2c}{b}\cosh(b\eta_1)}},
\label{eq:exact-integral_}
\end{equation}
where $\eta_1=\eta(0)$ and $q_1=\left.\frac{d\eta}{dx}\right|_{x=0}$. \eqref{eq:exact-integral_} defines an implicit relation between $x$ and $\eta(x)$. The integral is elliptic. It is often evaluated numerically. The parameters $\eta_1$ and $q_1$ are not directly specified, as discussed in Section~\ref{sec:model}. In practice, galvanostatic and potentiostatic setups more naturally provide either $(\eta_1,\eta_2)$ or $(q_1,q_2)$, with $\eta_2=\eta(W)$ and $q_2=\left.\frac{d\eta}{dx}\right|_{x=W}$. We therefore use the shooting method \cite{ascher1995numerical}. It determines $\eta_1$ when $(q_1,q_2)$ are known, or determines $q_1$ when $(\eta_1,\eta_2)$ are known.

\begin{figure}[h!]
    \centering
        \includegraphics[width=1\linewidth]{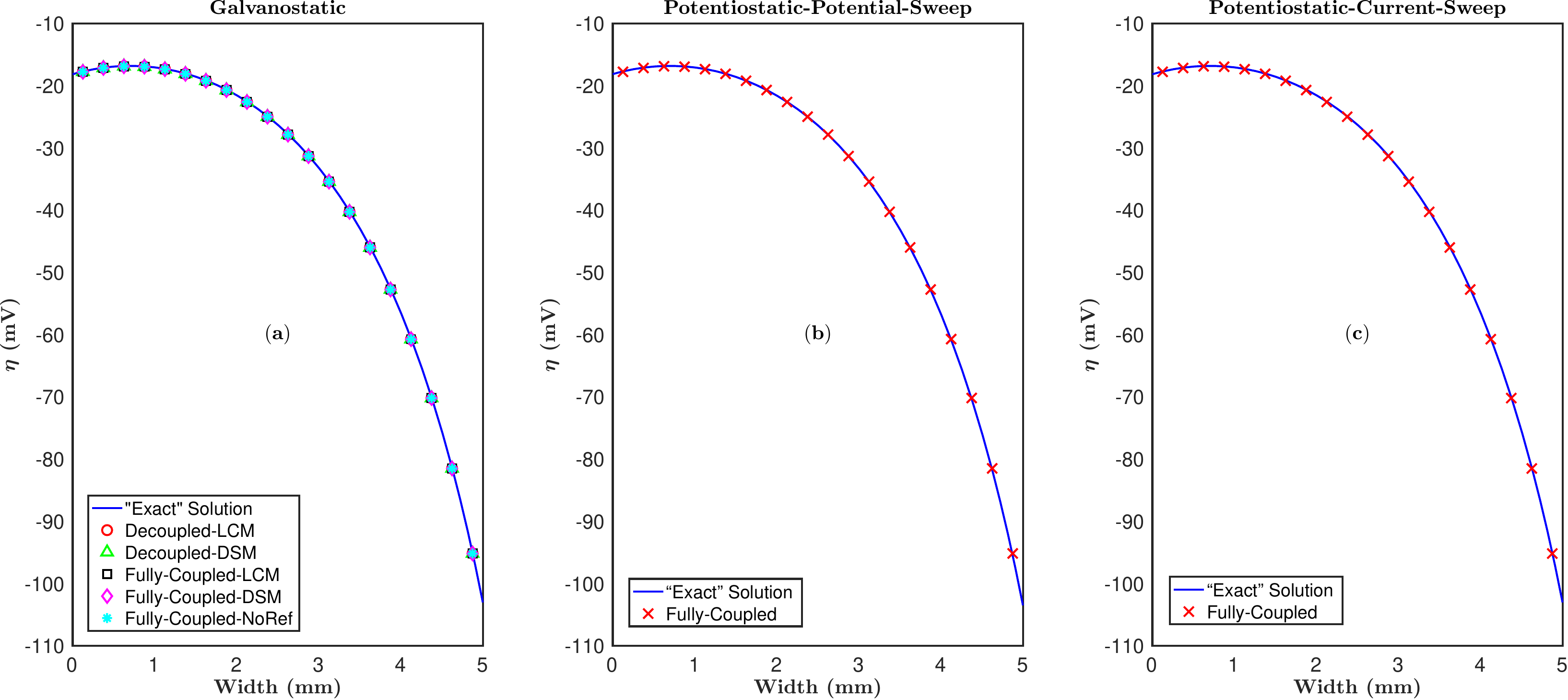}
        \caption{Numerical vs ``exact" solutions. (a): Galvanostatic; (b): Potentiostatic with potential sweep; (c): Potentiostatic with current sweep}
    \label{fig:exact-comparison}
\end{figure}

\begin{figure}[ht!]
    \centering
    \includegraphics[width=1\linewidth]{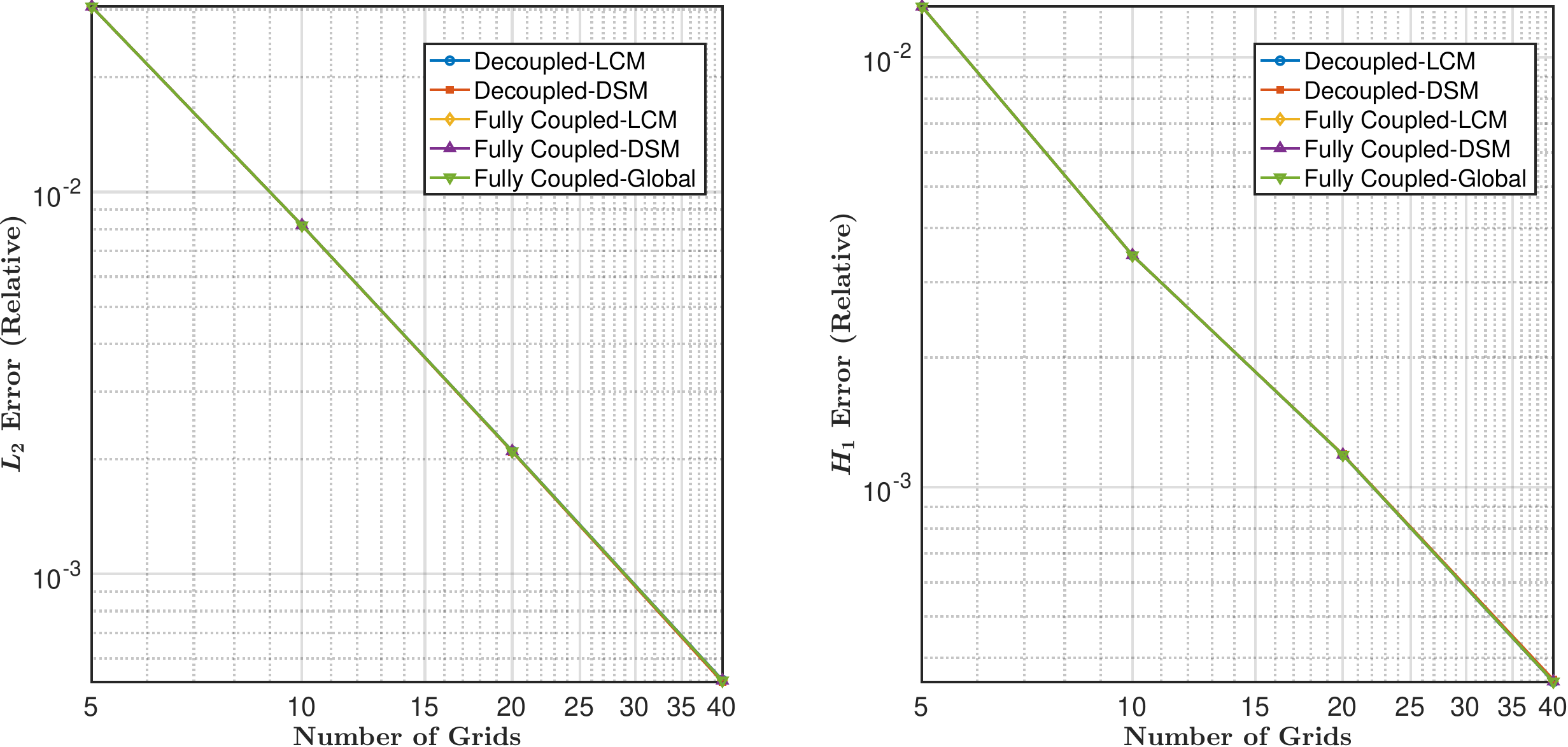}
    \caption{Error reduction ($L_2$, $H_1$)}
    \label{fig:error-reduction}
\end{figure}

For the galvanostatic case, \eqref{eq:BCe-g} and \eqref{eq:BCl-g} give $q_1=\frac{j_{\text{applied}}}{\sigma}$ and $q_2=-\frac{j_{\text{applied}}}{\kappa}$. For the potentiostatic case, we consider two setups. One sweeps boundary current at $x=W$, so $q_2$ is known. In this case, charge conservation \eqref{eq:dirichlet-replace} yields $q_1=-\frac{\kappa}{\sigma}q_2$. The other setup sweeps boundary potential at $x=W$. Then we lack $(\phi_l(0),\phi_e(W))$ to compute $(\eta_1,\eta_2)$. We need the numerical solver in Section~\ref{sec:methods-potentiostatic} to obtain these boundary values.
Figure~\ref{fig:exact-comparison} compares the numerical solutions with the corresponding ``exact'' solutions. Figure~\ref{fig:error-reduction} shows the $L_2$ and $H_1$ error reduction under mesh refinement. The numerical solutions match the ``exact'' profiles closely. The errors decrease quadratically with grid size. This confirms the accuracy and correctness of the proposed methods. 

\subsection{Objective Function for Finding the Electrolyte Reference Potential}
\begin{figure}[ht!]
    \centering
    \includegraphics[width=0.6\linewidth]{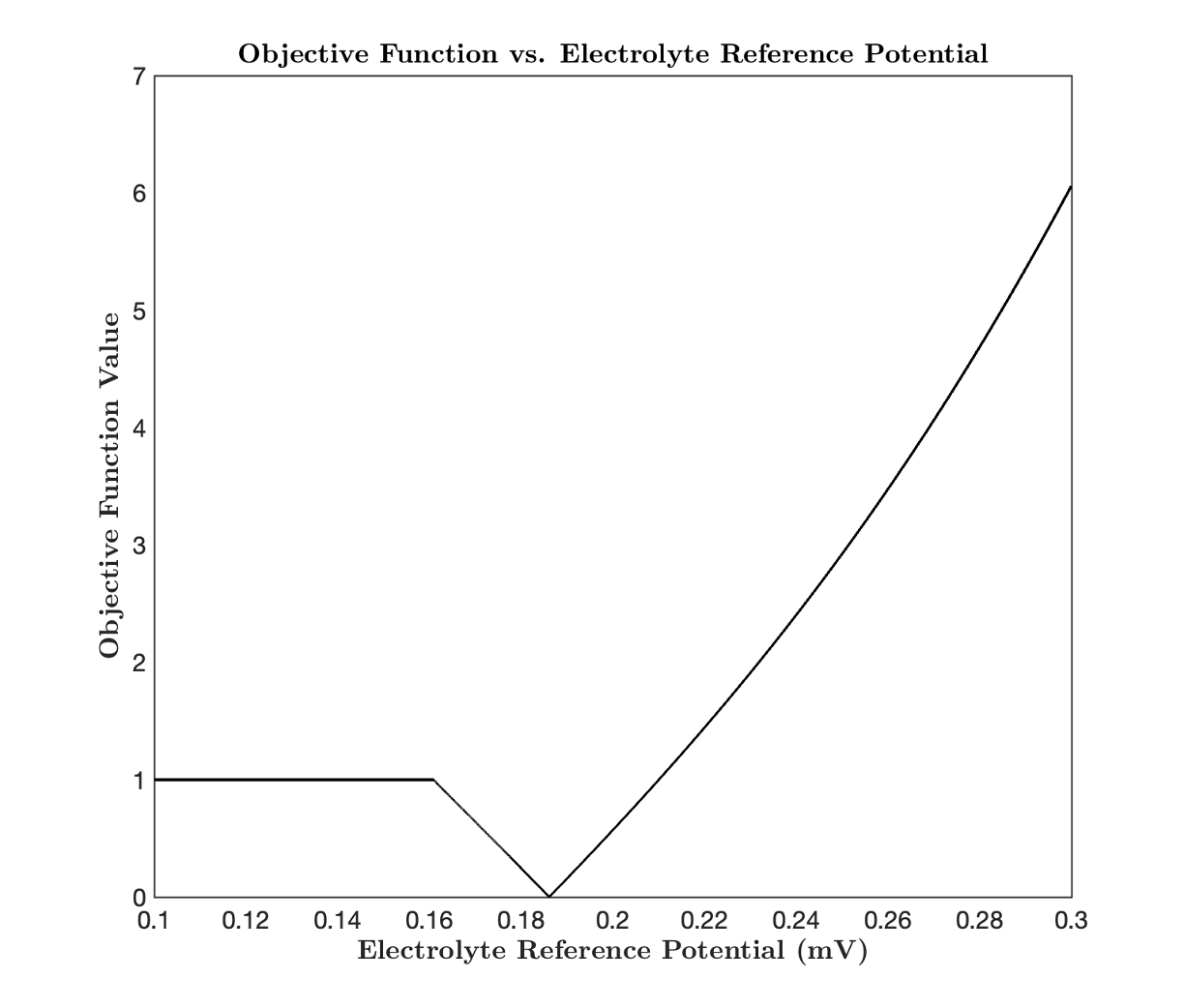}
    \caption{Objective function}
    \label{fig:obj}
\end{figure}

\vspace{-1.8ex}
As mentioned in Section~\ref{sec:methods-galvanostatic}, in a decoupled scheme the reference potential in $\Omega_l$ must be determined self-consistently during the solution process. We treat this as a search problem. The search is guided by the objective function in \eqref{eq:objective_function}, which enforces compatibility between the boundary flux and the sink/source term. No simple analytical expression is available, so we approximate its discrete realization numerically by sampling the parameter space of the electrolyte reference potential. Figure~\ref{fig:obj} shows the objective function for $j_{\text{applied}}=100~\mathrm{A/m^2}$. Extensive sampling over different $j_{\text{applied}}$ values indicates a consistent piecewise-like trend. The function is flat for smaller reference potentials on the left of the minimum. It rises sharply after the minimum, indicating strong sensitivity to the reference potential. Right-side deviations lead to amplified discrepancies. The ``V''-like shape near the minimum suggests that quadratic methods such as Newton's method may be unsuitable, since second-derivative information is unlikely to be helpful there. Gradient descent may also struggle near the tip and can oscillate without proper step-size adaptation. A carefully designed step-size adjustment strategy is beneficial to avoid such issues. Alternatively, derivative-free optimization methods, such as Golden-section Search \cite{press2007numerical}, may also perform well for this type of objective function.

\subsection{Homogeneous vs. Heterogeneous}
\begin{figure}[h!]
    \centering
    \includegraphics[width=.95\linewidth]{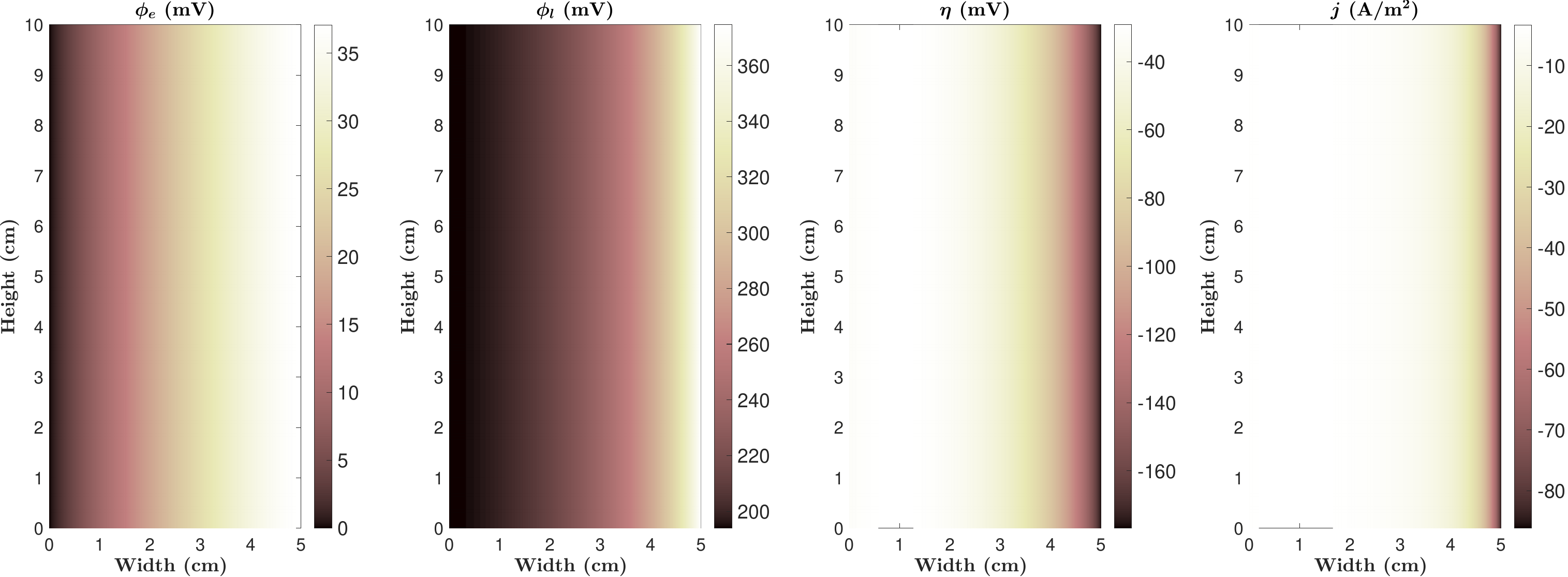}
    \caption{Potential (\ensuremath{\phi_e}, \ensuremath{\phi_l}, \ensuremath{\eta}) and local current density distributions (\ensuremath{j}) [Galvanostatic,  
    \ensuremath{j_{\text{applied}} = 1000 \,\text{A/m}^2}, homogeneous \ensuremath{\sigma} and \ensuremath{\kappa}, x-axis is scaled by 10 times]}
    \label{fig:homogeneous_50x50}
\end{figure}

\begin{figure}[h!]
    \centering
    \includegraphics[width=1\linewidth]{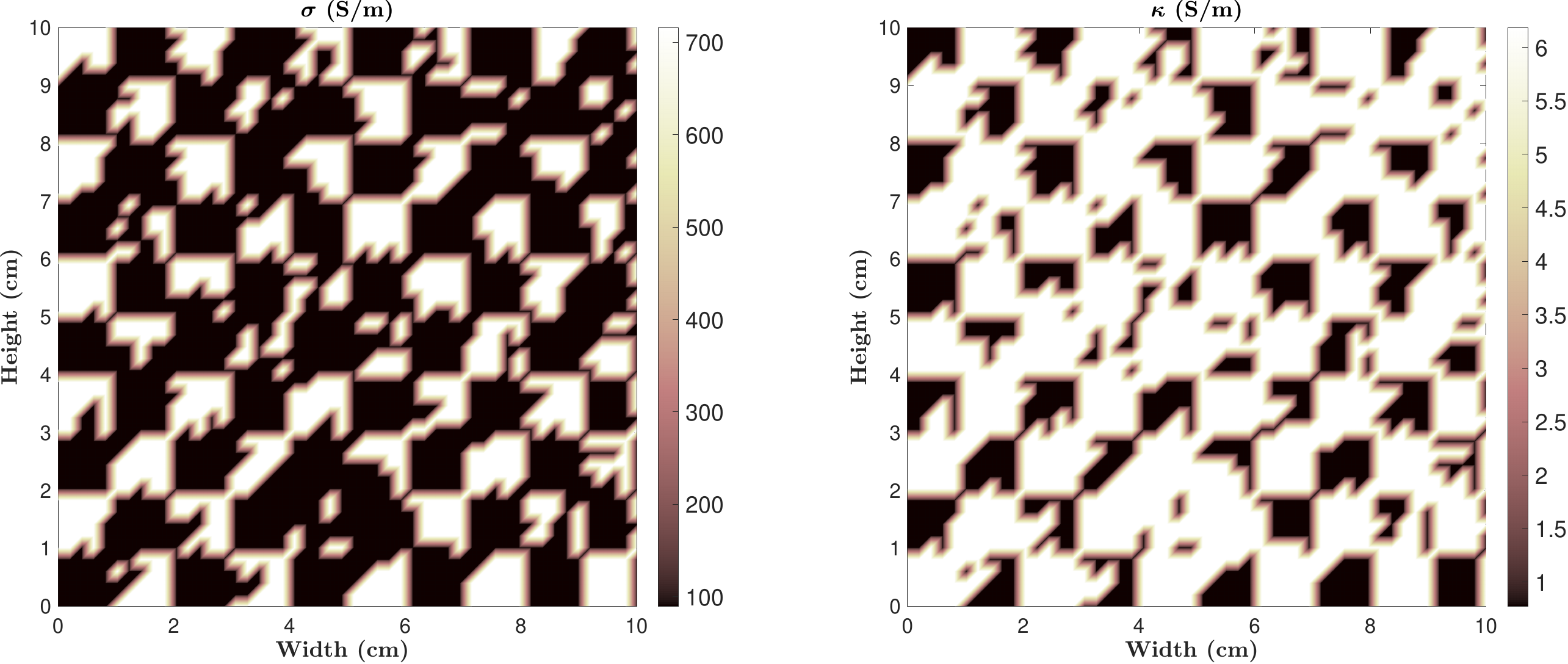}
    \caption{Bimodal porous electrode conductivity (\ensuremath{\sigma}) field and electrolyte conductivity (\ensuremath{\kappa}) field [x-axis is scaled by 20 times]}
    \label{fig:bimodal}
\end{figure}

\begin{figure}[h!]
    \centering
    \includegraphics[width=1\linewidth]{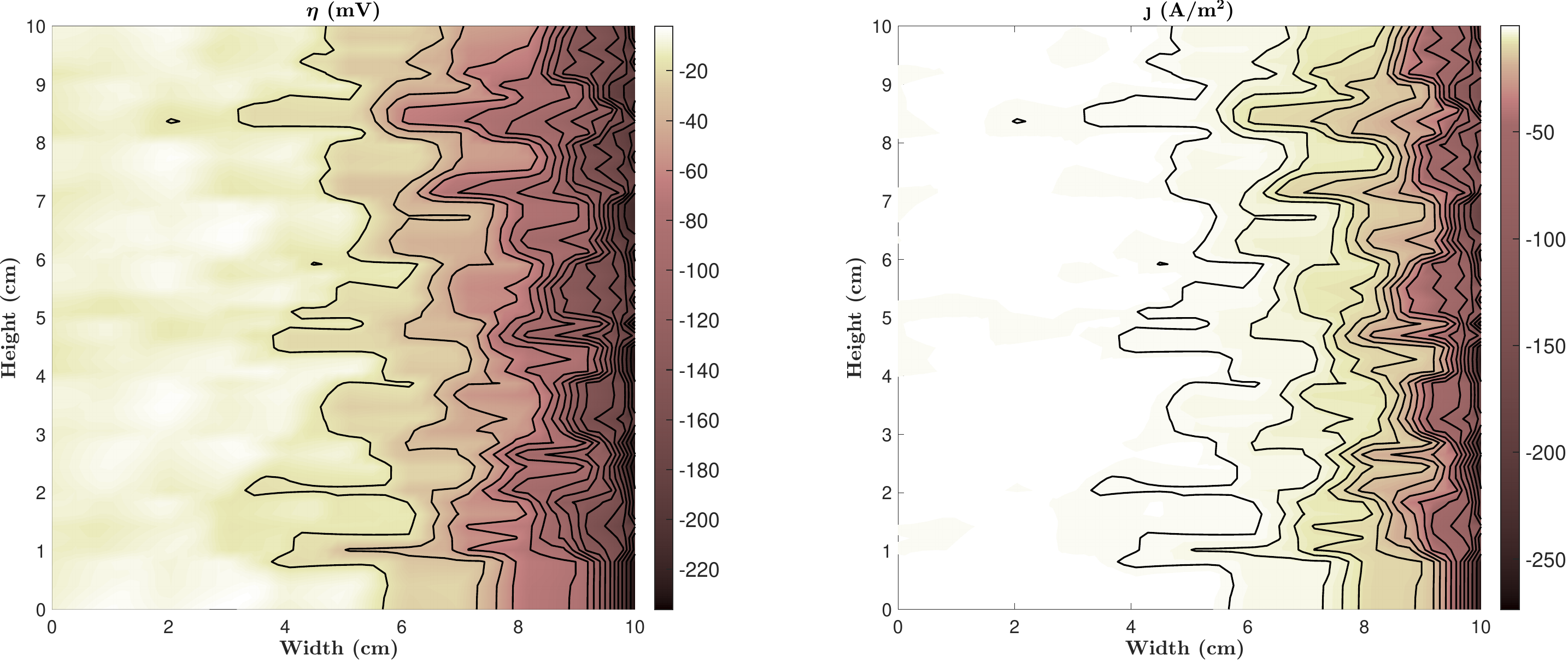}
    \caption{Overpotential (\ensuremath{\eta}) field and local current density (\ensuremath{j}) field of the bimodal conductivity case [\(j_\text{applied} = 1000 \, \text{A/m}^2\), x-axis is scaled by 20 times]}
    \label{fig:bimodal_eta_j}
\end{figure}

\begin{figure}[h!]
    \centering
    \includegraphics[width=1\linewidth]{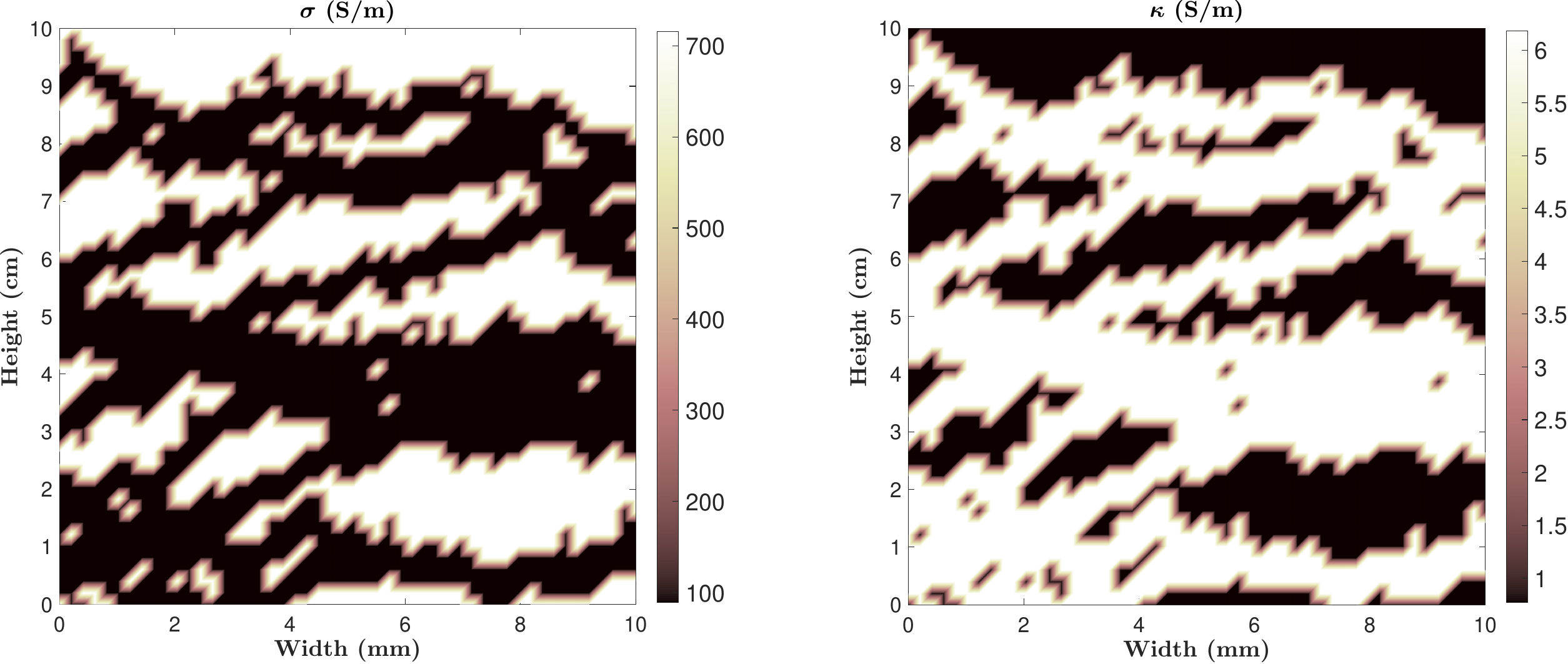}
    \caption{Channelized porous electrode conductivity (\ensuremath{\sigma}) field and electrolyte conductivity (\ensuremath{\kappa}) field [x-axis is scaled by 20]}
    \label{fig:channelized}
\end{figure}

\begin{figure}[h!]
    \centering
    \includegraphics[width=1\linewidth]{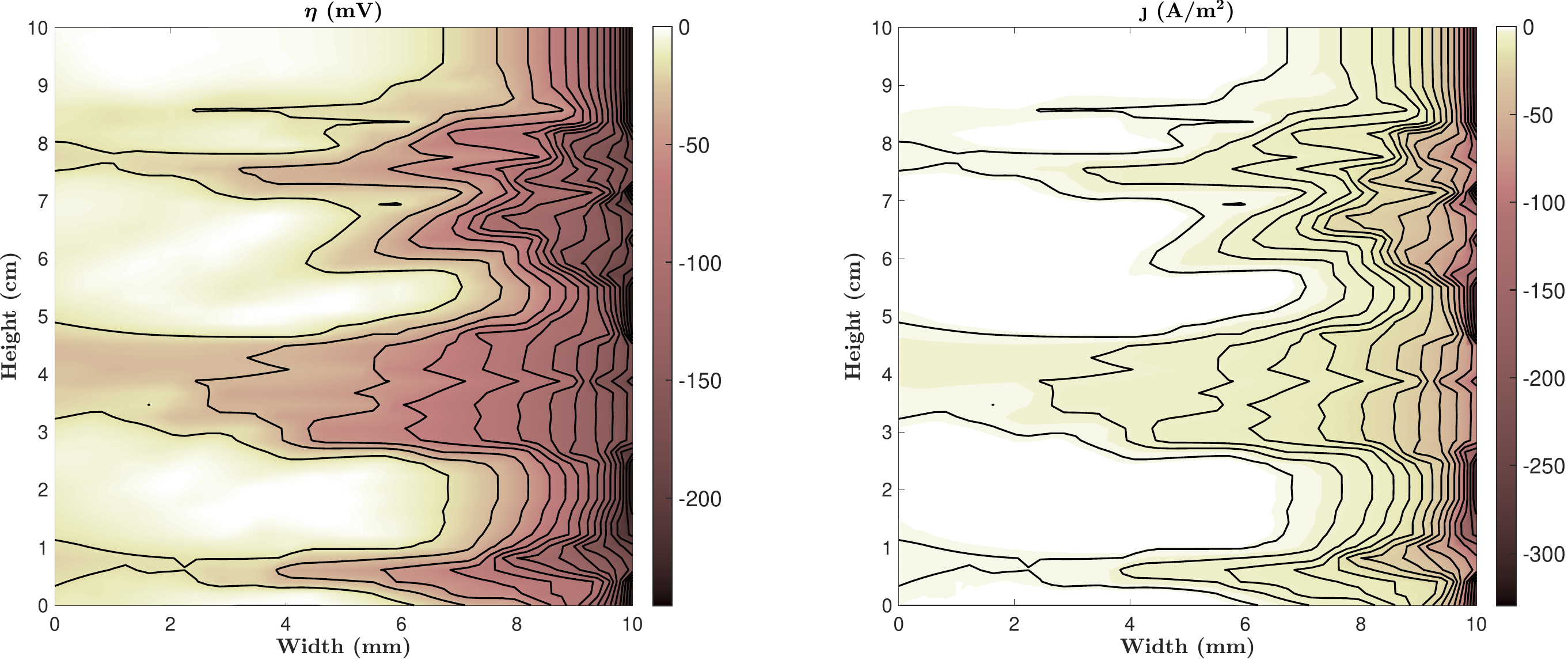}
    \caption{Overpotential (\ensuremath{\eta}) field and local current density (\ensuremath{j}) field of the channelized conductivity case [\(j_\text{applied} = 1000 \, \text{A/m}^2\), x-axis is scaled by 20]}
    \label{fig:channelized_eta_j}
\end{figure}

\vspace{-1.8ex}
The potential and local current density distributions for homogeneous and heterogeneous conductivity fields are illustrated in Figures~\ref{fig:homogeneous_50x50}, \ref{fig:bimodal_eta_j}, and \ref{fig:channelized_eta_j}. 
The corresponding bimodal and channlized conductivity fields are shown in Figures~\ref{fig:bimodal} and \ref{fig:channelized}. Since the results obtained from different methods are very close, we present only the galvanostatic results obtained using fully coupled LCM scheme for the case of $j_{\text{applied}}=1000~\mathrm{A/m^2}$. The computational grid is structured with dimension of $50\times50$.

In the homogeneous case, $\phi_e$, $\phi_l$, $\eta$, and $j$ vary smoothly. $\phi_l$ is much larger than $\phi_e$ because $\sigma$ is about 17 times larger than $\kappa$, leading to a larger potential drop in the electrolyte domain. $\eta$ is smallest near the current collector ($x=0$) and largest near the separator ($x=W$), confirming reaction dominance near the separator. The current density $j$ is strongly localized near the separator, driven by the electrolyte potential drop. The patterns match physical expectations and validate the solver.

For heterogeneous tests, we use randomized bimodal conductivity fields for $\sigma$ and $\kappa$ based on a bimodal porosity field with $\epsilon_{\text{low}}=0.2$ and $\epsilon_{\text{high}}=0.8$. Figure~\ref{fig:bimodal} shows one realization, and Figure~\ref{fig:bimodal_eta_j} shows the resulting $\eta$ and $j$ fields. Compared with the homogeneous case, the contours become irregular due to the bimodal structure of $\sigma$ and $\kappa$. Reaction localization remains near the separator (right boundary), where electrolyte-dominated transport yields higher overpotentials. The jagged features in $j$ indicate current crowding, highlighting the impact of heterogeneity and connectivity on charge transport efficiency.
Figure~\ref{fig:channelized_eta_j} shows the corresponding $\eta$ and $j$ distributions for the channelized conductivity case. The plots show localized regions of intensified overpotential and current density induced by channel connectivity. Together with the bimodal case, these results demonstrate that the proposed numerical methods can handle heterogeneous conductivity fields.

\subsection{Decoupled vs. Fully Coupled}
\begin{figure}[h!]
    \centering
    \includegraphics[width=1\linewidth]{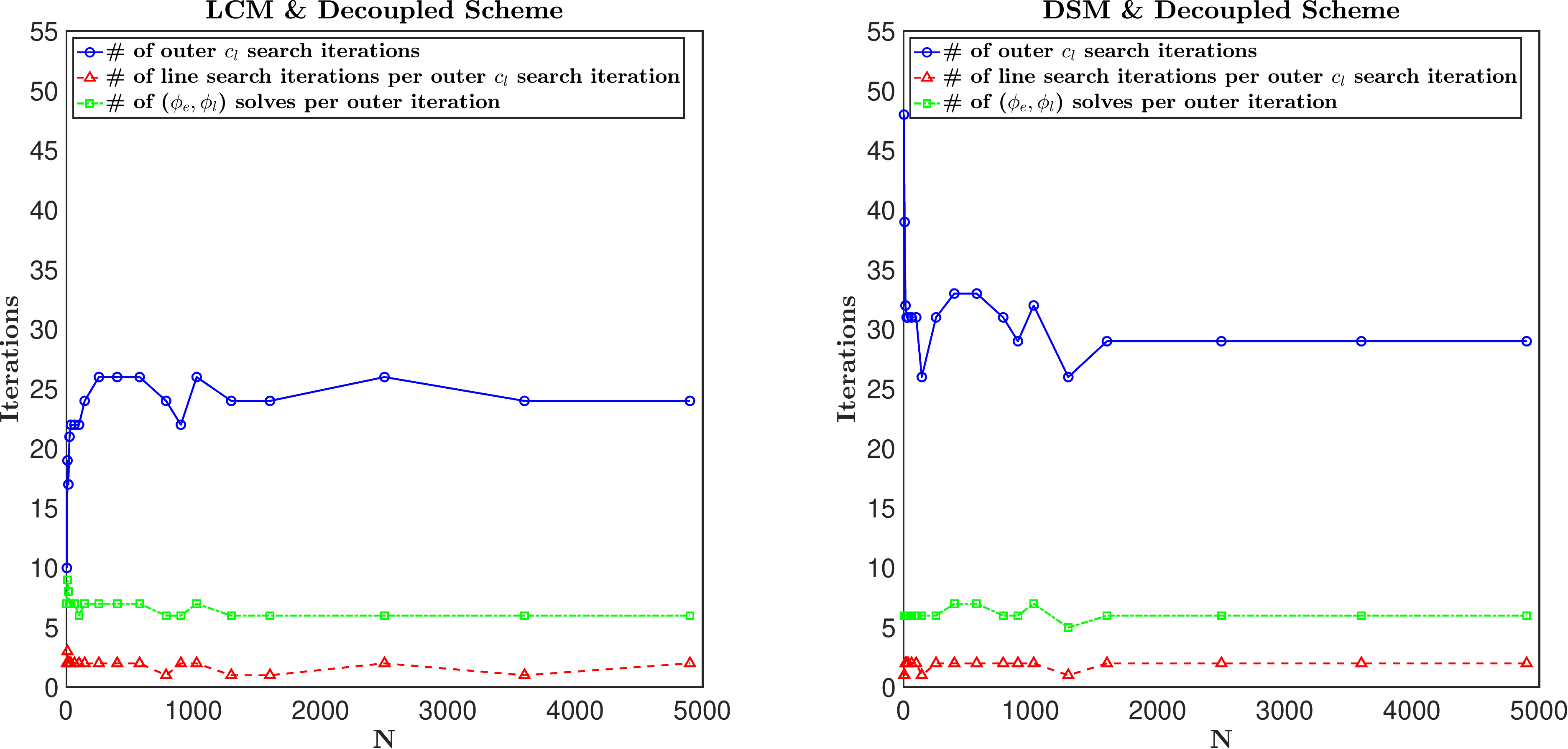}
    \caption{Iteration performance of the decoupled Scheme using LCM and DSM (\(j_\text{applied} = 500 \, \text{A/m}^2\), $N$ is the size of Jacobian Matrix. The global residual tolerance is \(1 \times 10^{-6}\))}
    \label{fig:iterations_vs_grids}
\end{figure}

\begin{figure}[h!]
    \centering
    \includegraphics[width=1\linewidth]{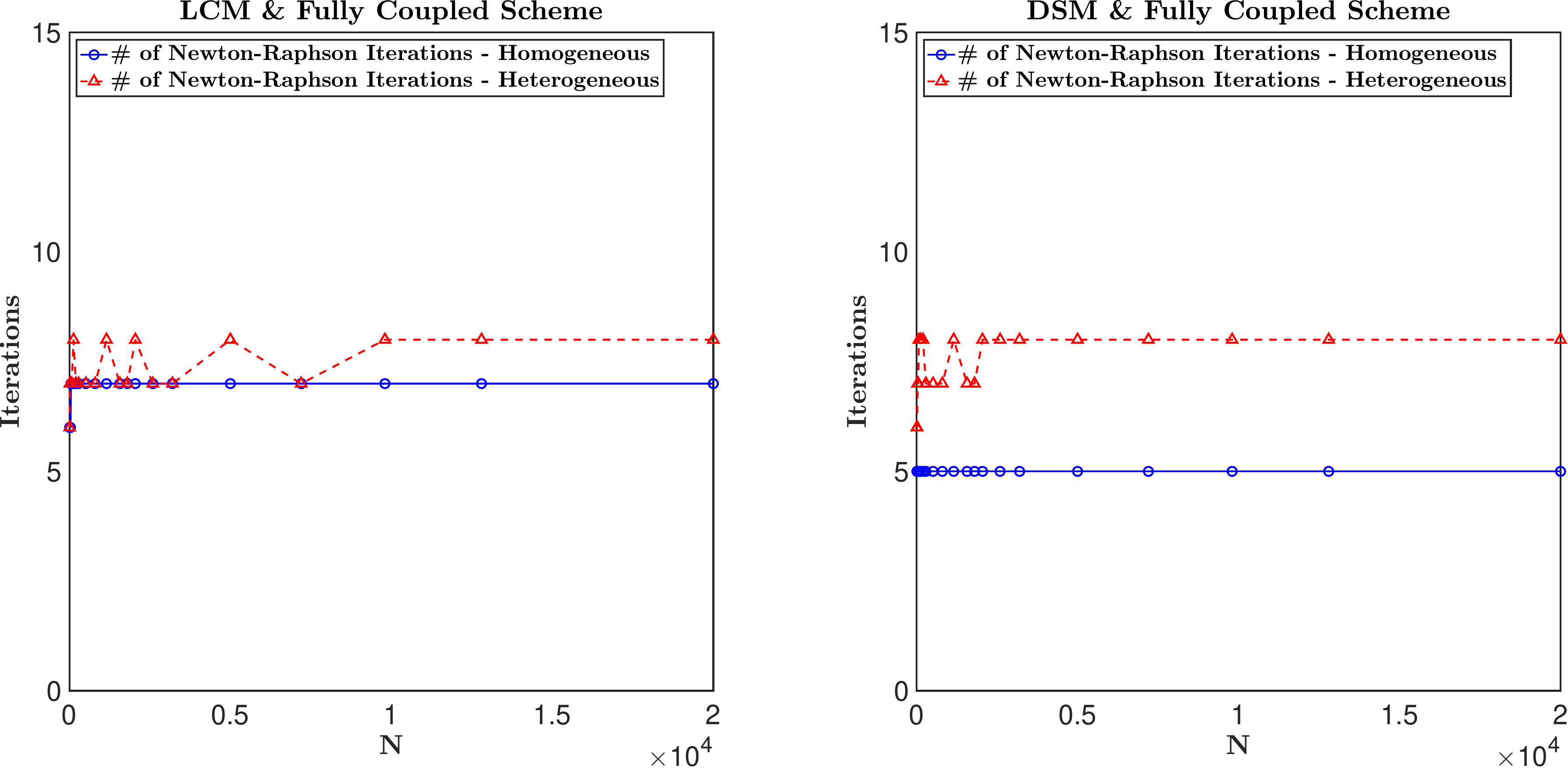}
    \caption{Iteration performance of the fully coupled scheme using LCM and DSM (\(j_\text{applied} = 500 \, \text{A/m}^2\), N is the size of Jacobian Matrix. The global residual tolerance is \(1 \times 10^{-6}\). The heterogeneity is the bimodal case)}
    \label{fig:iterations_vs_grids_NR}
\end{figure}

\vspace{-1.8ex}
As described in Section~3.3.1, the decoupled scheme requires an outer search for the electrolyte reference potential $c_l$ to maintain self-consistency with the coupled dynamics. Figure~\ref{fig:iterations_vs_grids} shows that refining the grid does not change the convergence behavior of the $c_l$ search. The search iteration count is approximately independent of $N$ for both LCM and DSM. In these runs, LCM typically uses about $20\%$ fewer search iterations than DSM. The drawback of the decoupled scheme is the added computational cost of this search. Each update of $c_l$ requires multiple line-search iterations, and the two Poisson equations are solved two times per line-search iteration \cite{bonnans2006numerical}. Each sequential solve also requires multiple nonlinear iterations. The linear systems solved in the decoupled scheme are smaller, with each system about half the size of the fully coupled system. However, the outer search repeats these solves many times and can introduce an overhead of up to hundreds of times the cost of the fully coupled approach. Therefore, despite its simpler implementation, the decoupled scheme is not practical due to its high computational expense.
In contrast, the fully coupled scheme solves the two Poisson equations only once. Figure~\ref{fig:iterations_vs_grids_NR} shows strong scalability for the fully coupled schemes with LCM and DSM. In these runs, DSM outperforms LCM for the homogeneous case and performs roughly the same as LCM for the heterogeneous case. 

It should be noted that, for the decoupled solution, a full nonlinear iteration is not necessary, since the outer search only needs an estimate of the search gradient. We therefore perform one nonlinear iteration per $(\phi_e,\phi_l)$ solve, and the final solution is identical to runs with full iterations. Nevertheless, as long as such searching is used, the fully coupled scheme remains more attractive.
The Jacobian-based iterations are sensitive to the initialization of $(\phi_e,\phi_l)$. The initial values affect the residual vector and the derivative terms in the Jacobian. These terms include highly nonlinear exponential functions and can become extremely large in magnitude, which may lead to divergence. For stability, we recommend initializing $\phi_e=\phi_l=0$. Alternatively, set $\phi_e=0$ and $\phi_l=-E_{eq}$.
LCM and DSM show similar convergence in our tests, suggesting that DSM may be preferable. LCM adds one internal constant-potential constraint, so $\mathrm{size}(J^L)\approx \mathrm{size}(J)$. However, it can produce a poorly scaled $J^L$, which may require special rescaling to reduce its condition number and can increase the number of linear solver iterations.

\subsection{LCM/DSM vs. GCM}
\begin{figure}[ht!]
    \centering
    \includegraphics[width=1\linewidth]{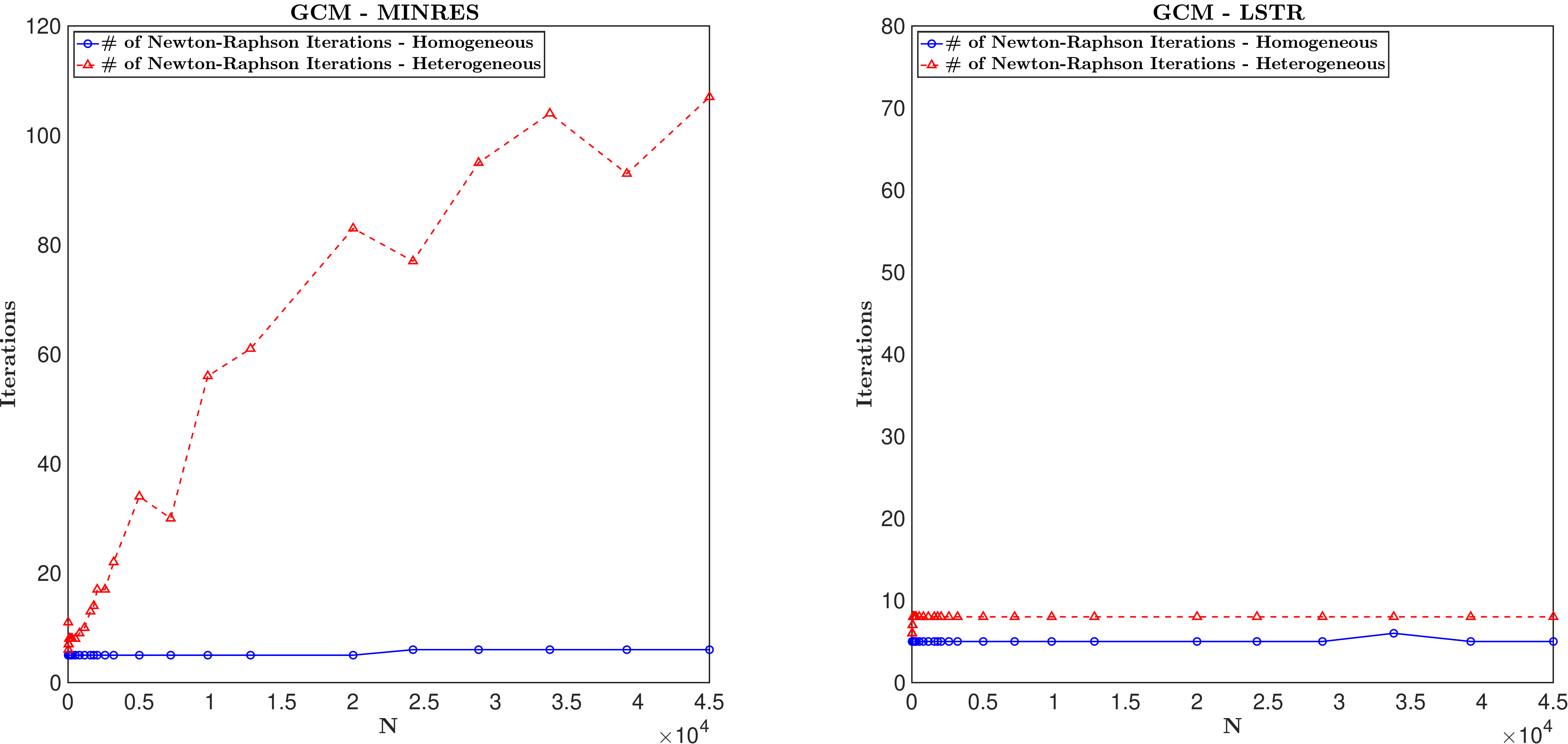}
    \caption{Iteration performance of the fully coupled scheme using MINRES and LSTR ($j_\text{applied}=500~\mathrm{A/m}^2$, $N$ is the size of the Jacobian matrix, global residual tolerance $1\times10^{-6}$, heterogeneity is the bimodal case).}
    \label{fig:iterations_vs_grids_global}
\end{figure}

\vspace{-1.8ex}
To assess the fully coupled scheme with GCM, we consider cases without imposing an explicit constant reference potential and report the convergence behavior for MINRES and LSTR in Figure~\ref{fig:iterations_vs_grids_global}. We do not include the Moore--Penrose pseudoinverse, because its $\mathcal{O}(N^3)$ cost can quickly exhaust computational resources. For the homogeneous case, MINRES and LSTR show identical performance, and they are comparable to the LCM/DSM results in Figure~\ref{fig:iterations_vs_grids_NR}. For the heterogeneous case, MINRES struggles to reach the prescribed tolerance $1\times10^{-6}$, although the outer nonlinear iterations still converge. In contrast, LSTR performs well for the heterogeneous case. Therefore, LSTR is preferable to MINRES. We notice that inexact linear solves can still support global nonlinear convergence, as long as the linear residual reaches about $10^{-2}$--$10^{-3}$. This motivates iterative solvers with moderate accuracy for large-scale problems. However, inexact linear solves can degrade nonlinear convergence from quadratic to linear and may introduce mild oscillations.

\begin{figure}[ht!]
    \centering
    \includegraphics[width=0.6\linewidth]{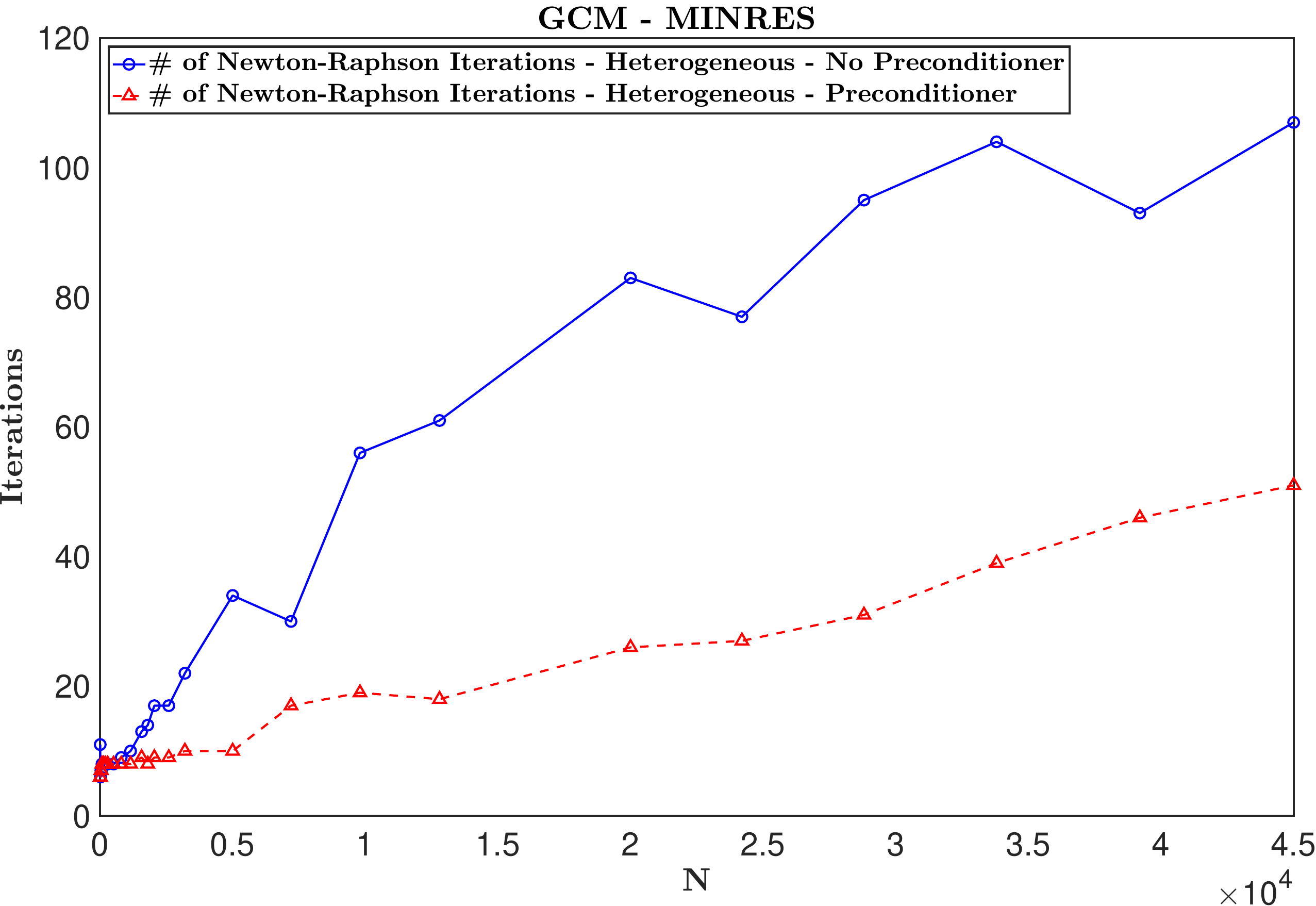}
    \caption{Iteration performance of the fully coupled scheme using MINRES with preconditioning (\(j_\text{applied} = 500 \, \text{A/m}^2\). \(N\) represents the size of the Jacobian matrix. The global residual tolerance is \(1 \times 10^{-6}\). The heterogeneity is the bimodal case)}
    \label{fig:iterations_vs_grids_prec}
\end{figure}

To improve MINRES, we apply preconditioning. The Jacobian $J$ is singular with $\mathrm{rank}=N-1$ and has a one-dimensional null space. We use the augmented system
\begin{equation}
    J_{\text{aug}} =
    \begin{bmatrix}
    J & \lambda I \\
    \lambda I & I
    \end{bmatrix}, \label{eq:J_aug}
\end{equation}
where $\lambda$ is a regularization parameter, typically in $10^{-6}$--$10^{-2}$. Equation~\eqref{eq:J_aug} shifts the singular values away from zero and reduces the condition number. The matrix $J_{\text{aug}}$ is twice the size of $J$, which can add considerable computational cost. To reduce this overhead, we use the Schur complement and obtain the reduced operator $J-\lambda^2 I$. Figure~\ref{fig:iterations_vs_grids_prec} shows that this preconditioning reduces the iteration count. The iteration count is still higher than LSTR. It also increases with $N$.

It should be noted that, for homogeneous conductivity fields, GCM with MINRES or LSTR is a viable approach, as it avoids the need to impose explicit potential references during the solve; the pair $(\phi_e, \phi_l)$ is then recovered via the postprocessing shift in Eq.~\eqref{eq:postprocessing}. For heterogeneous cases, the computational cost can be substantially higher than that of LCM/DSM when MINRES is used, and remains elevated with LSTR as well. Nevertheless, the key point is that galvanostatic solutions for $\phi_e - \phi_l$ or $\eta$ can be obtained without specifying any explicit potential reference.

\section{Conclusions} \label{sec:conclusions}
\vspace{-1.8ex}
This work systematically provides the numerical methods for solving the nonlinearly coupled Poisson equations in dual-continuum modeled porous-electrode systems. Under galvanostatic conditions, we show that the system is uniquely determined in terms of the potential difference $\phi_e-\phi_l$ (or $\eta$), while $\phi_e$ and $\phi_l$ share a free additive constant, which makes the discrete system underconstrained and singular. To address this, we present three strategies: LCM, DSM, and GCM. LCM is based on Dirichlet’s principle and a Lagrange-multiplier constraint derived from a unified energy functional for the coupled Poisson dynamics. DSM follows from Gauss’s theorem. LCM and DSM impose an explicit reference potential, whereas GCM solves for $\phi_e-\phi_l$ (or $\eta$) without specifying a reference.

We also present decoupled and fully coupled nonlinear solution schemes. The fully coupled scheme is more efficient and robust, because the decoupled scheme requires an outer search for a reference potential to maintain self-consistency. In homogeneous cases, LCM, DSM, and GCM yield comparable performance. In heterogeneous cases, GCM can struggle with MINRES, even with preconditioning, while GCM with LSTR remains stable and efficient. Numerical results validate the methods across multiple conductivity fields and operating conditions, including underconstrained singularities, nonlinear source/sink coupling, and heterogeneity. 
The provided numerical framework is not restricted to 2D setting. It can be extended to 3D and to more general multicontinuum advection--diffusion--reaction models that fully couple species transport with electrochemical kinetics. We view these extensions as future directions to further the endeavor in multiphysics porous-electrode modeling.

\section*{Data and Code Availability}
\vspace{-1.8ex} 
The code and data used in this study are available in the GitHub repository \url{https://github.com/harrywang1129/porous_electrode_solver}.

\bibliographystyle{unsrt}
\bibliography{references}

\appendix
\section*{Supplementary Information}
\renewcommand{\thesubsection}{\Alph{subsection}}
\renewcommand{\theequation}{\Alph{subsection}.\arabic{equation}}
\makeatletter
\@addtoreset{equation}{subsection}
\makeatother

\subsection{Solution Uniqueness of Galvanostatic Solution}
\label{sec:uniqueness}
\setcounter{equation}{0} 
\vspace{-1.8ex}
We prove that the following system admits a unique set of solution pairs \((\phi_e, \phi_l)\) up to an additive constant:
\setlength{\jot}{10pt} 
\begin{subequations} 
\begin{align}
\nabla \cdot (-\sigma \, \nabla \phi_e) &= -f(\eta) \quad \text{in  } \Omega_e, \label{eq:poisson-11_appendix}\\
\nabla \cdot (-\kappa \, \nabla \phi_l) &= \phantom{-}f(\eta) \quad \text{in  } \Omega_l,  \label{eq:poisson-22_appendix} \\
\sigma \frac{\partial \phi_e}{\partial n} &= g_e \quad \text{on  } \partial\Omega_e, \label{eq:BCe-g_appendix} \\
\kappa \frac{\partial \phi_l}{\partial n} &= g_l \quad \text{on  } \partial\Omega_l.
\label{eq:BCl-g_appendix}
\end{align}
\end{subequations}
where \( f = a \, \sinh(b \, \eta) \) and \(\eta = \phi_e - \phi_l - E_{eq}\). 

Let \((\phi_e^1, \phi_l^1)\) and \((\phi_e^2, \phi_l^2)\) be two pairs of solutions to \eqref{eq:poisson-11_appendix}--\eqref{eq:BCl-g_appendix}. We subtract \eqref{eq:poisson-11_appendix} and \eqref{eq:poisson-22_appendix} for \(\phi_{e/l}^2\) from \mwang{those equations} for \(\phi_{e/l}^1\) and define \(\xi_e := \phi_e^1 - \phi_e^2\), \(\xi_l := \phi_l^1 - \phi_l^2\), \(\eta^{1} := \phi_e^1 - \phi_l^1 -E_{eq}\), and \(\eta^{2} :=  \phi_e^2 - \phi_l^2-E_{eq}\). They satisfy the following equations:
\begin{equation*}
    \begin{split}
            \nabla\cdot(-\sigma \nabla \xi_e ) &= -f(\eta^1) + f(\eta^2),\\
            \nabla\cdot(-\kappa \nabla \xi_l ) &= f(\eta^1) - f(\eta^2).\\
    \end{split}
\end{equation*}
Here, the subscripts \((e/l)\) refer to the electrode and electrolyte, respectively, while the superscripts \((1/2)\) denote the first and second pairs of solutions.

We then multiply both equations by \((\xi_e, \xi_l)\) respectively and perform integration over \(\Omega_{e/l}\), we have:
\setlength{\jot}{10pt} 
\begin{subequations} 
\begin{align}
    \int_{\Omega_{e}} \xi_e \,\nabla \cdot (-\sigma \nabla \xi_e) \, dx &= - \int_{\Omega_{e}} \xi_e \,[f(\eta^1) - f(\eta^2)] \, dx,
    \label{eq:integral_xi_e} \\
    \int_{\Omega_{l}} \xi_l \,\nabla \cdot (-\kappa \nabla \xi_l) \, dx &= \int_{\Omega_{l}} \xi_l \, [f(\eta^1) - f(\eta^2)] \, dx.
    \label{eq:integral_xi_l}
\end{align}
\end{subequations}
We process the integral terms on the left-hand side by integration by parts and we can arrive at:
\setlength{\jot}{10pt} 
\begin{subequations} 
\begin{align*}
    \int_{\Omega_{e}} \xi_e \, \nabla \cdot (-\sigma \nabla \xi_e) \, dx 
    &= \int_{\Omega_{e}} \sigma |\nabla \xi_e|^2 \, dx
    - \int_{\partial \Omega_{e}} \xi_e \sigma \frac{\partial \xi_e}{\partial n} \, dS_e, \\
    \int_{\Omega_{l}} \xi_l \,\nabla \cdot (-\kappa \nabla \xi_l) \, dx 
    &= \int_{\Omega_{l}} \kappa |\nabla \xi_l|^2 \, dx
    - \int_{\partial \Omega_{l}} \xi_l \kappa \frac{\partial \xi_l}{\partial n} \, dS_l,
\end{align*}
\end{subequations}
Since \(\xi_e := \phi_e^1 - \phi_e^2\), \(\xi_l := \phi_l^1 - \phi_l^2\), we know \(\sigma \frac{\partial \xi_e}{\partial n} = 0\) \mwang{on $\partial\Omega_e$} and  \(\kappa \frac{\partial \xi_l}{\partial n} = 0\) \mwang{on $\partial\Omega_l$}. And thus, 
\setlength{\jot}{10pt} 
\begin{subequations} 
\begin{align}
    \int_{\Omega_{e}} \sigma |\nabla \xi_e|^2 \, dx &= - \int_{\Omega_{e}} \xi_e \, [f(\eta^1) - f(\eta^2)] \, dx,
    \label{eq:integral_xi_e_} \\
    \int_{\Omega_{l}} \kappa |\nabla \xi_l|^2 \, dx &= \int_{\Omega_{l}} \xi_l \, [f(\eta^1) - f(\eta^2)] \, dx.
    \label{eq:integral_xi_l_}
\end{align}
\end{subequations}
By summing \eqref{eq:integral_xi_e_} and \eqref{eq:integral_xi_l_} together, we obtain:
\begin{equation}
    \begin{split}
       0\leq\int_{\Omega_{e}} \sigma |\nabla \xi_e|^2 +\int_{\Omega_{l}} \kappa |\nabla \xi_l|^2 &= - \int_{\Omega} [f(\eta^1) - f(\eta^2)](\xi_e -\xi_l)\ dx\\
       & = - \int_{\Omega} [f(\eta^1) - f(\eta^2)](\eta^1 -\eta^2)\ dx.
    \end{split} \label{eq:sum}
\end{equation}
Noticing that 
\begin{equation*}
    f'(\eta) = a b \,\cosh(b \, \eta) >0,
\end{equation*}
we conclude that \( f(\eta) \) is a strictly monotonically increasing function of \(\eta\). This implies that for any given \( x \),
\begin{equation*}
    \begin{split}
    \text{if } \eta^1(x)> \eta^2(x) &\implies f\Big(\eta^1(x)\Big) > f\Big(\eta^2(x)\Big),\\
    \text{if } \eta^1(x)< \eta^2(x) &\implies f\Big(\eta^1(x)\Big) < f\Big(\eta^2(x)\Big),\\
    &\Downarrow\\
    \Big[f\big(\eta^1(x)\big) - f\big(\eta^2(x)\big)\Big]&\Big(\eta^1(x) -\eta^2(x)\Big)\geq 0, \quad \forall x\in \Omega.
    \end{split}
\end{equation*}

Integrating over the entire domain, we obtain
\begin{equation*}
    - \int_{\Omega} \Big[f(\eta^1) - f(\eta^2)\Big] (\eta^1 -\eta^2)\ dx \leq 0.
\end{equation*}
Since the integral is non-negative, the only possibility for equality to hold is when the integrand itself is zero \mwang{almost everywhere} in \( \Omega \), leading to
\begin{equation}
    f(\eta^1) - f(\eta^2) = 0 \quad \text{a.e. in } \Omega.
\end{equation}
Since \( f(\eta) \) is strictly monotonic, this directly implies that
\begin{equation*}
    \eta^1(x) = \eta^2(x) \quad \text{a.e. in } \Omega.
\end{equation*}
Recalling that \( \eta = \phi_e - \phi_l - E_{\text{eq}} \), we immediately obtain
\begin{equation*}
    \xi_e(x) = \xi_l(x) \quad \text{a.e. in } \Omega.
\end{equation*}
Now, looking back at \eqref{eq:sum}, we notice
\begin{equation*}
    \int_{\Omega} \sigma |\nabla \xi_e|^2 \,dx + \int_{\Omega} \kappa |\nabla \xi_l|^2 \,dx  = 0,
\end{equation*}
we conclude that
\begin{equation*}
    \nabla \xi_e = 0, \quad \nabla \xi_l = 0 \quad \text{a.e. in } \Omega.
\end{equation*}
 This implies that both \(\xi_e\) and \(\xi_l\) must be spatially constant in \(\Omega\), i.e.,
\begin{equation*}
     \xi_e(x) = \xi_l(x) = C,\quad\text{a.e. in } \Omega
\end{equation*}
for some constant \(C \in \mathbb{R}\). 

Thus, any solution pair \((\phi_e, \phi_l)\) to the system must be of the form:
\begin{equation*}
    \{(\phi_e,\phi_l) \mid \phi_e = \phi_e^* + C ,\quad \phi_l = \phi_l^* + C, \quad C\in \mathbb{R} \}.
\end{equation*}
where the asterisk (\(*\)) denotes a particular pair of solutions. That is, \(\phi_e^*\) and \(\phi_l^*\) represent a specific set of solutions to the system, and any other solutions can be obtained by adding a constant \(C\).

This completes the proof that the solution pair to the coupled nonlinear Neumann Poisson system is unique up to an additive constant \(C\).

\subsection{Dirichlet's Principle and Lagrange Multiplier} 
\label{sec: Lagrange}
\setcounter{equation}{0} 
\vspace{-1.8ex}
We consider here a coupled Poisson system under all-Neumann boundary conditions in the following general form:
\setlength{\jot}{7pt} 
\begin{subequations}
\begin{align}
    \nabla \cdot (-\sigma \, \nabla \phi_1) &= -f(\phi_1 - \phi_2)  \quad \text{in } \Omega,  \label{eq:coupled1} \\ 
    \nabla \cdot (-\kappa \, \nabla \phi_2) &= f(\phi_1 - \phi_2) \quad \text{in } \Omega, \label{eq:coupled2} \\ 
    \sigma \, \frac{\partial \phi_1}{\partial n} &= g_1  \quad \text{on } \partial \Omega,  \label{eq:coupled-neumann-bc-1} \\  
    \kappa \, \frac{\partial \phi_2}{\partial n} &= g_2 \quad \text{on } \partial \Omega. \label{eq:coupled-neumann-bc-2}
\end{align}
\end{subequations}
where \(f(\phi_1 - \phi_2) = a \, \sinh{[b \, (\phi_1 - \phi_2-E_{eq})]}\), and \(a\) and \(b\) are some constants.

Based on Dirichlet's principle, We try to write a unified energy functional for this system as
\begin{equation}
    \mathcal{E}(\mwang{\phi_1}, \mwang{\phi_2}) = \frac{1}{2} \int_\Omega \sigma \, |\nabla \phi_1|^2 + \frac{1}{2} \int_\Omega \kappa \, |\nabla \phi_2|^2 \mwang{-} \int_\Omega V(\phi_1,\phi_2) - \int_{\partial \Omega} g_1 \phi_1 - \int_{\partial \Omega} g_2 \phi_2  \label{eq:unified-energy}
\end{equation}
where:
\begin{itemize}
    \item \(\frac{1}{2} \int_\Omega \sigma \, |\nabla \phi_1|^2\) and \(\frac{1}{2} \int_\Omega \kappa \, |\nabla \phi_2|^2\) represents the respective energy associated with the gradients of \(\phi_1\) and \(\phi_2\), penalizing variations in \(\phi_1\) and \(\phi_2\) within \(\Omega\);
    \item \(V(\phi_1,\phi_2)\) is a potential energy term accounting for the interaction between the source/sink term (\(-f, f\)) and (\(\phi_1, \phi_2\)) within \(\Omega\);
    \item \(\int_{\partial \Omega} g_1 \phi_1\) and  \(\int_{\partial \Omega} g_2 \phi_2\) incorporates the contributions of the Neumann boundary fluxes \(g_1\) and \(g_2\) to the energy of the coupled system.     
\end{itemize}
For \(V(\phi_1,\phi_2)\), we have that its partial derivative with respect to \(\phi_1\) yields \(-f(\phi_1-\phi_2)\) and with respect to \(\phi_2\) yields \(f(\phi_1-\phi_2)\). In a matrix form, we can express as:
\begin{equation}
\renewcommand{\arraystretch}{2.6} 
\setlength{\arraycolsep}{4pt}    
\text{\normalsize$
\begin{bmatrix}
    \frac{\partial V(\phi_1,\phi_2)}{\partial \phi_1} \\
    \frac{\partial V(\phi_1,\phi_2)}{\partial \phi_2} 
\end{bmatrix}
$}
=
\text{\normalsize$
\begin{bmatrix}
    -f(\phi_1-\phi_2) \\ f(\phi_1-\phi_2)
\end{bmatrix}.
$}
\label{eq:potential-function}
\end{equation}
\textbf{Note}: For the case of (3a)--(3b), a natural choice for such potential function is
\begin{equation}
    V(\phi_e,\phi_l) = -\frac{a}{b} \, \cosh(b \, \eta)
\end{equation}
with \(\eta = \phi_e - \phi_l - E_{eq}\).

Next, we perform variational analysis to \eqref{eq:unified-energy} by introducing a small perturbation to \(\phi_1\) such that:
\begin{equation}
    \phi_1^{\epsilon_1} = \phi_1 + \epsilon_1 v_1, \label{eq:perturbation-phi1} 
\end{equation}
where \(\epsilon_1\) is an arbitrary small variable, and \(v_1\) is an arbitrary smooth function. 

Plug \eqref{eq:perturbation-phi1} in, we have:
\begin{equation}
\begin{split}
    \mathcal{E}(\phi_1^{\epsilon_1}, \phi_2) &= \frac{1}{2} \int_\Omega \sigma \, |\nabla (\phi_1 + \epsilon_1 v_1)|^2 
    + \frac{1}{2} \int_\Omega \kappa \, |\nabla \phi_2|^2  \yw{-} \int_\Omega V[(\phi_1 + \epsilon_1 v_1),\phi_2] \\
    &- \int_{\partial \Omega} g_1 (\phi_1 + \epsilon_1 v_1) 
    - \int_{\partial \Omega} g_2 \phi_2  
\end{split}
\label{eq:unified-energy_phi1}
\end{equation}
Here \(|\nabla (\phi_1 + \epsilon_1 v_1)|^2 = |\nabla \phi_1|^2 + 2 \epsilon_1 \nabla \phi_1\cdot \nabla v_1 + \epsilon_1^2 |\nabla v_1|^2\). 

Differentiate \eqref{eq:unified-energy_phi1} with respect to \(\epsilon_1\) and evaluate at \(\epsilon_1 = 0\), then: 
\setlength{\jot}{10pt} 
\begin{subequations}
\begin{align}
    \frac{\partial}{\partial \epsilon_1} \left( \frac{1}{2} \int_\Omega  \sigma \, |\nabla (\phi_1 + \epsilon_1 v_1)|^2 \right) \bigg|_{\epsilon_1=0} &=  \int_\Omega \sigma \, \nabla \phi_1 \cdot \nabla v_1 \\
    &= \int_\Omega \nabla \cdot (- \sigma \, \nabla \phi_1) \, v_1 + \int_{\partial \Omega} \sigma \, \frac{\partial \phi_1}{\partial n} v_1 , \label{eq:sqauregradient}\\
    \frac{\partial}{\partial \epsilon_1} \left( \int_\Omega V[(\phi_1 + \epsilon_1 v_1),\phi_2]\right) \bigg|_{\epsilon_1=0} &= \int_\Omega \mwang{-}f(\phi_1-\phi_2) \, v_1, \label{eq:f_}\\
    \frac{\partial}{\partial \epsilon_1} \left( \int_{\partial \Omega} g_1 (\phi_1 + \epsilon_1 v_1) \right) \bigg|_{\epsilon_1=0} &= \int_{\partial \Omega} g_1 v_1.
\end{align}
\end{subequations}
Similarly, if we perform the variation with respect to \(\phi_2\), we have:
\setlength{\jot}{10pt} 
\begin{subequations}
\begin{align}
    \frac{\partial}{\partial \epsilon_2} \left( \frac{1}{2} \int_\Omega  \kappa \, |\nabla (\phi_2 + \epsilon_2 v_2)|^2 \right) \bigg|_{\epsilon_2=0} &=  \int_\Omega \kappa \, \nabla \phi_2 \cdot \nabla v_2 \\
    &= \int_\Omega \nabla \cdot (- \kappa \, \nabla \phi_2) \, v_2 + \int_{\partial \Omega} \kappa \, \frac{\partial \phi_2}{\partial n} v_2 , \label{eq:sqauregradient_}\\
    \frac{\partial}{\partial \epsilon_2} \left( \int_\Omega V[\yw{\phi_1,}(\phi_2 + \epsilon_2 v_2)]\right) \bigg|_{\epsilon_2=0} &= \int_\Omega f(\phi_1-\phi_2) \, v_2, \label{eq:f__}\\
    \frac{\partial}{\partial \epsilon_2} \left( \int_{\partial \Omega} g_2 (\phi_2 + \epsilon_2 v_2) \right) \bigg|_{\epsilon_2=0} &= \int_{\partial \Omega} g_2 v_2.
\end{align}
\end{subequations}
Thus, we have:
\yw{
\setlength{\jot}{10pt}
\begin{equation} \label{eq:E_epsilon}
    \left.\frac{\partial \mathcal{E}(\phi_1^{\epsilon_1}, \phi_2^{\epsilon_2})}{\partial (\epsilon_1, \epsilon_2)}\right|_{\epsilon_1=0,\, \epsilon_2=0}
    =
    \begin{bmatrix}
        \int_\Omega \nabla \cdot (- \sigma \, \nabla \phi_1) \, v_1 + \int_\Omega f(\phi_1-\phi_2) \, v_1 - \int_{\partial \Omega} (g_1 - \sigma \, \frac{\partial \phi_1}{\partial n})v_1  \\
          \int_\Omega \nabla \cdot (- \kappa \, \nabla \phi_2) \, v_2 - \int_\Omega f(\phi_1-\phi_2) \, v_2 - \int_{\partial \Omega} (g_2 - \kappa \, \frac{\partial \phi_2}{\partial n})v_2
    \end{bmatrix}   
\end{equation}  
}
To minimize the energy functional \(\mathcal{E}(\phi_1, \phi_2)\), we set \eqref{eq:E_epsilon} to zero such that: 
\yw{
\begin{equation} \label{eq:unified_energy=0}
    \frac{\partial \mathcal{E}(\phi_1^{\epsilon_1}, \phi_2^{\epsilon_2})}{\partial (\epsilon_1, \epsilon_2)} \bigg|_{\epsilon_1=0, \, \epsilon_2=0} 
    = 
    \begin{bmatrix}
    0  \\
    0
    \end{bmatrix}
\end{equation} 
}
Since \(v_1\) and \(v_2\) are arbitrary, \eqref{eq:unified_energy=0} hold if and only if:
\setlength{\jot}{7pt} 
\begin{subequations}
\begin{align}
    \nabla \cdot (-\sigma \, \nabla \phi_1) &= -f(\phi_1 - \phi_2)  \quad \text{in } \Omega,   \\ 
    \nabla \cdot (-\kappa \, \nabla \phi_2) &= f(\phi_1 - \phi_2) \quad \text{in } \Omega, \\ 
    \sigma \, \frac{\partial \phi_1}{\partial n} &= g_1  \quad \text{on } \partial \Omega,   \\  
    \kappa \, \frac{\partial u_2}{\partial n} &= g_2 \quad \text{on } \partial \Omega. 
\end{align}
\end{subequations}
We can see that the above obtained Euler--Lagrange equations exactly recover \eqref{eq:coupled1}--\eqref{eq:coupled2}, and honor the original Neumann boundary conditions as in \eqref{eq:coupled-neumann-bc-1}--\eqref{eq:coupled-neumann-bc-2}. By constructing a unified energy functional in this manner, even when the coefficients \(\sigma\) and \(\kappa\) differ, the variational formulation recovers the original coupled system in its entirety.

\textbf{Remark:} If the nonlinear term lacks antisymmetry, then no potential function \( V(\phi_1, \phi_2) \) can reproduce \eqref{eq:potential-function}, provided the reaction current density depends on \( \phi_1 - \phi_2 \)
 In other words, the system cannot be derived from a single unified energy functional when the coefficients of the nonlinear terms in the two equations are not balanced. This violates the principle of ``conservativity" or ``self-adjointness." Consequently, if such coupled system is governed by an asymmetric nonlinear source/sink term, it is inherently impossible to construct a single energy functional to describe it.

Now, based on the unified energy functional \eqref{eq:unified-energy}, let us add some local Lagrange constraints to the coupled system. Suppose we want to impose constraints on \(\phi_1\) at certain specified positions \( \Gamma_1 \) within \(\Omega\), expressed as  
\begin{equation}
    \phi_1|_{\Gamma_1} = \vec{C}_1 \quad \text{at } \Gamma_1.
\end{equation}
We can use the Lagrange multiplier method to introduce a set of multipliers \( \lambda_{i} \) (which are additional global unknowns but only act at the ``constraint point"), and construct the following augmented energy functional:
\begin{equation}  
\mathcal{L}(\phi_1, \phi_2, \lambda) = \mathcal{E}(\phi_1, \phi_2) + \sum_{i=1}^m \lambda_{i} [\phi_1(x_i) - c_i] \label{eq:extended-functional-coupled}
\end{equation}
\(x_i \in \Omega\) are the discrete points to enforce constraint, given by the set \(\Gamma_1: =\{x_1, x_2, \dots, x_m\}\). 

To derive the extended Euler-Lagrange equations of \eqref{eq:extended-functional-coupled}, we take variation of \( \mathcal{L} \) with respect to \( \phi_1 \), \( \phi_2 \), and \(\lambda\). Utilizing the respective Euler-Lagrange equations of \(\mathcal{E}(u_1, u_2)\), we can obtain the following extended Lagrange system:
\setlength{\jot}{10pt} 
\begin{subequations}
\begin{align}
    \nabla \cdot (-\sigma \, \nabla \phi_1) + f + \sum_{i=1}^{m}\lambda_i \, \delta (x-x_i)& = 0  \quad \text{in } \Omega  \label{eq:extended-lagrange-1} \\
    \nabla \cdot (-\kappa \, \nabla \phi_2) - f  &= 0  \quad \text{in } \Omega  \label{eq:extended-lagrange-2} \\
   \phi_1(x_i) - c_i &= 0 \quad \text{for } \mwang{x_i\in}\Gamma_1   \label{eq:extended-lagrange-4}   
\end{align}
\end{subequations}
with Neumann conditions given in \eqref{eq:coupled-neumann-bc-1}--\eqref{eq:coupled-neumann-bc-2}.

Here, \(\delta(x-x_i)\) is Dirac delta function, meaning \(\lambda_i\) effects the system only at \(x_i\).

\subsection{Line Search Method for Reference Potential Matching} 
\label{sec:line-search}
\setcounter{equation}{0} 
\vspace{-1.8ex}
To refine the reference potential \(c_l\) for the electrolyte domain, we employ a gradient descent-based line search methods within an optimization framework. This iterative approach aims to minimize the objective function \(J\), as defined in (18), ensuring that the solution adheres to charge conservation. Gradient descent updates the reference potential by moving in the direction of the negative gradient of the objective function. The gradient is computed using a finite-difference approximation:
\begin{equation}
    \nabla J(c_l) \approx \frac{J(c_l + \delta) - J(c_l)}{\delta}
\end{equation}
Alternatively, central finite differences can be used for higher accuracy, given by:
\begin{equation}
    \nabla J(c_l) \approx \frac{J(c_l + \delta) - J(c_l - \delta)}{2\delta}
\end{equation}
where \(\delta\) is a small perturbation. The iterative update rule is:
\begin{equation}
    c_l^{(k+1)} = c_l^{(k)} - \beta \nabla J(c_l^{(k)})
\end{equation}
where \(\beta > 0\) is the step size. To ensure convergence, the step size is adjusted using the Armijo condition:
\begin{equation}
    J(c_l^{(k+1)}) \leq J(c_l^{(k)}) + c \beta \nabla J(c_l^{(k)})(-\nabla J(c_l^{(k)}))
\end{equation}
where \(c \in (0, 1)\) is a small constant, and \(\beta\) is initialized to a default value and iteratively adjusted as follows:
1. Start with an initial \(\beta > 0\).
2. Reduce \(\beta\) by a factor of \(\rho \in (0, 1)\) (e.g., \(\rho = 0.5\)) until the Armijo condition is met.
3. Use the updated \(\beta\) for the next iteration.

\subsection{Exact Solution Expression (27)} 
\label{sec:exact-solution}
\setcounter{equation}{0} 
\vspace{-1.8ex}
We use the First Integral Method to find the exact solution of (27) in 1D under galvanostatic and potentiostatic conditions.

\textbf{Step 1: First Integral}

The given PDE is:
\begin{equation}
    \frac{d^2 \eta}{dx^2} = a' \sinh(b\eta).
\end{equation}
Multiply both sides by \(\frac{d\eta}{dx}\) to prepare for integration:
\begin{equation}
    \frac{d^2 \eta}{dx^2} \frac{d\eta}{dx} = a' \sinh(b\eta) \frac{d\eta}{dx}
\end{equation}

Rewriting the left-hand side as a total derivative:
\begin{equation}
    \frac{d}{dx} \left( \frac{1}{2} \left(\frac{d\eta}{dx}\right)^2 \right) = \frac{d}{dx} \left( \frac{a'}{b} \cosh(b\eta) \right)
\end{equation}
Integrate both sides once:
\begin{equation}
    \frac{1}{2} \left(\frac{d\eta}{dx}\right)^2 = \frac{a'}{b} \cosh(b\eta) + C \label{eq:exact}
\end{equation}
where \(C\) is the integration constant.

Then, we use the Neumann boundary conditions to determine \(C\).

\textbf{Step 2: Apply Boundary Conditions}

At \(x = 0\):
\begin{equation}
    \frac{1}{2} q_1^2 = \frac{a'}{b} \cosh(b\eta_0) + C
\end{equation}
At \(x = W\):
\begin{equation}
    \frac{1}{2} q_2^2 = \frac{a'}{b} \cosh(b\eta_W) + C
\end{equation}
where $q_1 = \frac{d\eta}{dx}(0), q_2 = \frac{d\eta}{dx}(W), \eta_0 = \eta(0)$ and $\eta_W = \eta(W).$
We can use either boundary equation to find the expression for \(C\), such that:
\begin{equation}
    C = \frac{1}{2} q_1^2 - \frac{a'}{b} \cosh(b\eta_0)
\end{equation}
\textbf{Step 3: Solve for \(x\)}

Rearrange:
\begin{equation}
    \frac{d\eta}{dx} = \sqrt{\frac{2a'}{b} \cosh(b\eta) + q_1^2 - \frac{2a'}{b} \cosh(b\eta_0)}.
\end{equation}
Separating variables:
\begin{equation}
    dx = \frac{d\eta}{\sqrt{\frac{2a'}{b} \cosh(b\eta) + q_1^2 - \frac{2a'}{b} \cosh(b\eta_0)}}. \label{eq:exact-integral}
\end{equation}
Integrate:
\begin{equation}
    x = \int_{\eta_0}^{\eta} \frac{d\eta}{\sqrt{\frac{2a'}{b} \cosh(b\eta) + q_1^2 - \frac{2a'}{b} \cosh(b\eta_0)}}.
\end{equation}
This defines an implicit relation between \(x\) and \(\eta(x)\), which needs to be evaluated numerically as we end up with an integral problem. We can then apply numerical integration or approximate using elliptic functions.

\end{document}